\documentclass[12pt]{amsart}
\usepackage{amssymb, latexsym, amsmath, amscd, oldlfont}
\usepackage[dvips]{graphicx}
\usepackage{psfrag}

\advance\hoffset by -0.5 truecm

\def\e{\'e}

\def\V{\widehat{V}}
\def\U{\widehat{U}}

\def\H{{\cal H}}
\def\HH{{\widehat{H}}}
\def\K{\widehat{K}}

\def\ppi^n{\widetilde{\pi}}
\def\al{{\alpha}}

\def\ep{{\epsilon}}
\def\ga{{\gamma}}

\def\ov{\overline}
\def\om{{\omega}}
\def\Om{{\Omega}}

\def\te{{\theta}}

\def\lg{\langle}
\def\rg{\rangle}

\def\vr{\varphi}

\def\pr{\prime}
\def\wt{\widetilde}
\def\wh{\widehat}

\def\ppi{\widetilde{\pi}}

\def\lra{\longrightarrow}

\def\sm{\setminus}
\def\Sig{\Sigma}
\def\hra{\hookrightarrow}
\def\HZ{\operatorname{HZ}}

\def\C{{\mathbb C}}
\def\R{{\mathbb R}}
\def\N{{\mathbb N}}
\def\Z{{\mathbb Z}}

\newtheorem{theorem}{Theorem}[section]
\newtheorem{proposition}{Proposition}[section]
\newtheorem{lemma}{Lemma}[section]
\newtheorem{corollary}{Corollary}[section]
\theoremstyle{definition}
\newtheorem{definition}{Definition}[section]

\def\thebibliography#1{\section*{References}\list
 {[\arabic{enumi}]}{\settowidth\labelwidth{[#1]}\leftmargin\labelwidth
 \advance\leftmargin\labelsep
 \usecounter{enumi}}
 \def\newblock{\hskip .11em plus .33em minus .07em}
 \sloppy\clubpenalty4000\widowpenalty4000
 \sfcode`\.=1000\relax}

\newenvironment{remark}{\par\medskip\noindent{\bf 
Remark.}}{\hfill\par\medskip}
\newenvironment{remarks}{\par\medskip\noindent{\bf 
Remarks.}}{\hfill\par\medskip}
\newenvironment{question}{\par\medskip\noindent{\bf 
Question:}}{\hfill\par\medskip}

\author{Leonardo Macarini}
\title[Hofer-Zehnder capacity and Hamiltonian circle actions]{Hofer-Zehnder
capacity\\ and Hamiltonian circle actions}
\address{Instituto de Matem\'atica Pura e Aplicada - IMPA\\
         Estrada Dona Castorina, 110 - Jardim Bot\^anico\\
         22460-320 Rio de Janeiro RJ\\
         Brazil}
\email{leonardo@impa.br}         
\date{current version: April 2003}
\thanks{This work was partially supported by FAPERJ and CNPq-PROFIX, Brazil.}

\begin{document}

\begin{abstract}
We introduce the Hofer-Zehnder $G$-semicapacity $c_{\HZ}^G(M,\om)$ of a
symplectic manifold $(M,\om)$ (or $G$-sensitive Hofer-Zehnder capacity) with respect to a subgroup $G \subset \pi_1(M)$
($c_{\HZ}(M,\om) \leq c^G_{\HZ}(M,\om)$) and prove that given a geometrically bounded symplectic manifold $(M,\om)$ and an open subset $N \subset M$ admitting a Hamiltonian free circle
action with order greater than two then $N$ has bounded Hofer-Zehnder
$G$-semicapacity, where $G \subset \pi_1(N)$ is the subgroup generated by the
orbits of the action, provided that the index of rationality of $(M,\om)$ is
sufficiently great (for instance, if $[\om]|_{\pi_2(M)}=0$).

We give several applications of this result. Using P. Biran's decomposition
theorem, we prove the following: let $(M^{2n},\Om)$ be a closed K\"ahler
manifold ($n>2$) with $[\Om] \in H^2(M,\Z)$ and $\Sig$ a complex hypersurface
representing the Poincar\'e dual of $k[\Om]$, for some $k \in \N$. Suppose
either that $\Om$ vanishes on $\pi_2(M)$ or that $k>2$. Then there exists a
decomposition of $M$ into an open dense subset $E$ such that $E\sm\Sig$ has
finite Hofer-Zehnder $G$-semicapacity and an isotropic CW-complex, where $G
\subset \pi_1(E\sm\Sig)$ is the subgroup generated by the obvious circle action
on the normal bundle of $\Sigma$. Moreover, we prove that if $(M,\Sig)$ is
subcritical then $M\sm\Sig$ has finite Hofer-Zehnder $G$-semicapacity.

We also show that given a hyperbolic surface $M$ and $TM$ endowed with the
twisted symplectic form $\om_0 + \pi^*\Om$, where $\Om$ is the K\"ahler form on
$M$, then the Hofer-Zehnder $G$-semicapacity of the domain $U_k$ bounded by the
hypersurface of kinetic energy $k$ minus the zero section $M_0$ is finite if
$k < 1/2$, where $G \subset \pi_1(U_k)$ is the subgroup generated by the
fibers of $SM$.

Finally, we will consider the problem of the existence of periodic orbits on
prescribed energy levels for magnetic flows. We prove that given any weakly
exact magnetic field $\Om$ on any compact Riemannian manifold $M$ then there
exists a sequence of contractible periodic orbits of energy arbitrarily small,
extending a previous result of L. Polterovich.

\end{abstract}

\maketitle

\section{Introduction}

The problem of the existence of periodic orbits for Hamiltonian systems has a
very far and rich history in the development of Mathematics and Physics. In
fact, some of the major contributions on this question include the celebrated
Poincar\'e-Birkhoff theorem on the existence of periodic orbits for area
preserving twist maps on the annulus and the Morse-Lyusternik-Schnirelman 
theory on the existence of closed geodesics on compact Riemannian manifolds.

Recently, a special class of symplectic invariants was introduced by H. Hofer
and E. Zehnder \cite{HZ} as a link between rigidity phenomena of symplectic
mappings and the existence of periodic orbits provided by classical variational
methods. More precisely, they introduced an axiomatic definition of symplectic
capacities in the following way:

\begin{definition}
Consider the class of all symplectic manifolds $\{(M,\om)\}$ of fixed dimension
equal to 2n. A symplectic capacity is a map $(M,\om) \mapsto c(M,\om)$ that
associates to $(M,\om)$ a non negative number or $\infty$ and that satisfies the
following axioms:
\begin{itemize}
\item {\bf Monotonicity:} if there exists a symplectic embedding $\vr: (M,\om)
\to (N,\tau)$ then
$$ c(M,\om) \leq c(N,\tau); $$
\item {\bf Conformality:} $c(M,\al\om) =
|\al|c(M,\om)\ \forall\al\in\R\setminus\{0\};$
\item {\bf Non-triviality:} $c(B^{2n}(r),\om_0) = c(Z^{2n}(r),\om_0) = \pi
r^2,$
\end{itemize}
where $B^{2n}(r) = \{(x,y)\in\R^{2n}; \sum_{i=1}^n x_i^2+y_i^2\leq r\}$ is the
ball of radius $r$ and $Z^{2n}(r) = \{(x,y)\in\R^{2n}; x_1^2+y_1^2\leq r\}$ is
the cylinder over the $(x_1,y_1)$-plane.
\end{definition}

A distinguished symplectic capacity was introduced by Hofer and Zehnder and is
directly related to the question of the existence of periodic orbits on
prescribed energy levels of Hamiltonian flows:

\begin{definition}
\label{defcHZ}
Given a symplectic manifold $(M,\om)$, define the Hofer-Zehnder capacity of $M$
by
$$ c_{\HZ}(M,\om) = \sup_{H \in \H_a(M,\om)}\ \max H, $$
where $\H_a(M,\om)$ is the set of {\it admissible Hamiltonians} $H$ defined on
$M$, that is,
\begin{itemize}
\item $H \in \H(M,\om) \subset C^\infty(M,\R)$, where $\H(M,\om)$ is the set
of {\it pre-admissible Hamiltonians} defined on $M$, that is, $H$ satisfies the
following properties: $0 \leq H \leq m(H):=\max H$, there exist an open set 
$V \subset M$ such that $H|_V \equiv 0$ and a compact set 
$K \subset M\setminus\partial M$ satisfying $H|_{M\setminus K}\equiv m(H)$;
\item every non-constant periodic orbit of $X_H$ has period greater than 1.
\end{itemize}
\end{definition}

Intuitively, this capacity measure the sufficient oscillation for a
pre-admissible Hamiltonian to have fast periodic orbits, that is, of period 
less than one. It is easy to see that it satisfies the monotonicity and
conformality properties. Remarkably, by making use of classical variational
methods, Hofer and Zehnder \cite{HZ} showed that it also satisfies the
non-triviality axiom, giving, in particular, a dynamical new proof of Gromov's
nonsqueezing theorem \cite{Gro}.

Moreover, this symplectic invariant has important consequences on the question
of the existence of periodic orbits on prescribed energy levels of Hamiltonian 
systems. In fact, it is easy to prove that if a symplectic manifold $(M,\om)$
has {\it bounded Hofer-Zehnder capacity}, that is, if $c_{\HZ}(U,\om) < \infty$
for every open subset $U \subset M$ with compact closure, then given any
Hamiltonian $H: M \to \R$ with compact energy levels, there exists a dense
subset $\Sigma \subset H(M)$ such that for every $c \in \Sigma$ the energy
hypersurface $H^{-1}(c)$ has a periodic solution.

This result was improved by M. Struwe in \cite{Str} showing that if the energy
hypersurfaces bound a compact submanifold then $\Sigma$ has total Lebesgue
measure. This condition is not necessary as was proved recently in \cite{MS}. Note that as an immediate consequence we have a proof of the
Weinstein  conjecture \cite{We} for manifolds with bounded Hofer-Zehnder
capacity, since we can construct in a neighborhood of a contact hypersurface
$S$ a distinguished parametrized family of hypersurfaces with equivalent flows.

Unfortunately, it is, in general, a difficult task to determine when a
symplectic manifold has bounded Hofer-Zehnder capacity. For $(\R^{2n},\om_0)$, 
where $\om_0$ denotes the canonical symplectic form, it follows by
Hofer-Zehnder  theorem \cite{HZ}. For surfaces, it was proved by K. Siburg
\cite{Si}, that the  Hofer-Zehnder capacity of an open subset $U$ coincides
with the area of $U$.  However, the situation is dramatically different for
higher dimensional manifolds. Actually, Hofer and Viterbo showed in \cite{HV}
that 
$$ c_{\HZ}(M \times D^2(r),\om \oplus \om_0) = \pi r^2, $$ 
where $(M,\om)$ is any closed symplectic manifold and $D^2(r)$ is the two
dimensional disk of radius $r>0$, with $r$ satisfying the condition $\pi r^2 \leq
m(M,\om)$, where $m(M,\om)$ is the {\it index of rationality} of
$(M,\om)$, that is, 
$$ m(M,\om) = \inf\bigg\{\int_{[u]}\om;\ \ [u] \in \pi_2(M) \text{
satisfies } \int_{[u]}\om > 0\bigg\}. $$
If the above set is empty, we define $m(M,\om) = \infty$. Recall that a
2-form $\om$ on $M$ is {\it weakly exact} if it is closed and $m(M,\om) = \infty$.

This result was extended by G. Lu \cite{Lu1,Lu2} to noncompact symplectic
manifolds  satisfying certain conditions at the infinity. More precisely, he
showed that
$$ c_{\HZ}(M \times D^2(r),\om \oplus \om_0) \leq \pi r^2, $$
provided that $(M,\om)$ is geometrically bounded and $\pi r^2 \leq m(M,\om)$.

\begin{definition}
A symplectic manifold $(M,\om)$ is called {\it geometrically bounded} if there exists an almost complex structure $J$ on $M$ and a Riemannian metric $g$ such that
\begin{itemize}
\item $J$ is uniformly $\om$-tame, that is, there exist positive constants $c_1$ and $c_2$ such that
$$ \om(v,Jv) \geq c_1\|v\|^2 \text{ and } |\om(v,w)| \leq c_2\|v\|\|w\| $$
for all tangent vectors $v$ and $w$ to $M$;
\item the sectional curvature of $g$ is bounded from above and the injectivity radius is bounded away from zero.
\end{itemize}
\end{definition}

The closed symplectic manifolds are clearly geometrically bounded; a product of two geometrically bounded
symplectic manifolds is also such a manifold. Moreover, it was proved in
\cite{Lu1} that the cotangent bundle of a compact manifold endowed with any
twisted symplectic form is geometrically bounded (a proof of this fact can also be found in
\cite{CGK}). A {\it twisted symplectic structure} on a cotangent bundle $T^*M$
is a symplectic form given by $\om_0 + \pi^*\Om$, where $\om_0$ is the
canonical symplectic form, $\pi: T^*M \to M$ is the canonical projection and
$\Om$ is a closed 2-form on $M$.

Thus, using the monotonicity property, we can conclude that $(T^*(N\times
S^1),\om_0)$ has bounded Hofer-Zehnder capacity, where $N$ is any compact
manifold  and $\om_0$ is the canonical symplectic form (because we can embed
symplectically any bounded subset of $T^*M \times T^*S^1$ into $T^*M \times \R^2$). Actually, by the
same reason, we can prove that $(M \times T^*S^1,\om\oplus\om_0)$ has
bounded  Hofer-Zehnder capacity for any geometrically bounded symplectic manifold $M$. It should
be noted that it is necessary to consider product symplectic forms $\om
\oplus \om_0$. In fact, there is an example due to Zehnder \cite{Zeh} 
of a symplectic form on $T^3 \times \R \simeq T^2 \times
T^*S^1$ with unbounded Hofer-Zehnder capacity. Note that any symplectic form on
$T^2$ is weakly exact.

Other known example of cotangent bundle with bounded Hofer-Zehnder capacity
is $T^*T^n$ endowed with any twisted symplectic form $\om_0 +
\pi^*\Om$, where $\Om$ is a closed 2-form on $T^n$ and $\pi:T^*T^n \to T^n$ is
the canonical projection. It was proved by J. Jiang \cite{Ji} for
$\Om\equiv 0$ and by V. Ginzburg and  E. Kerman \cite{GK} for the general case.

We will consider here a refinement of the Hofer-Zehnder capacity by
considering periodic orbits whose homotopy class is contained in a given
subgroup $G$ of $\pi_1(M)$. More precisely, we define the Hofer-Zehnder 
semicapacity as follows; it can be compared with the
$\pi_1$-sensitive Hofer-Zehnder capacity defined by M. Schwarz \cite{Sch} and
Lu \cite{Lu1,Lu2} and the relative symplectic capacity introduced by P. Biran,
L. Polterovich and D. Salamon \cite{BPS}:

\begin{definition}
Given a symplectic manifold $(M,\om)$ and a subgroup $G \subset \pi_1(M)$
define the {\it Hofer-Zehnder G-semicapacity} of $M$ (or {\it $G$-sensitive Hofer-Zehnder capacity}) by
$$ c_{\HZ}^G(M,\om) = \sup_{H \in \H_a^G(M,\om)}\ \max H, $$
where $\H_a^G(M,\om)$ is the set of {\it $G$-admissible Hamiltonians} defined
on $M$, that is,
\begin{itemize}
\item $H \in \H(M,\om)$, that is, $H$ is pre-admissible (see definition
\ref{defcHZ});
\item every nonconstant periodic orbit of $X_H$ whose homotopy class belongs to
$G$ has period greater than 1.
\end{itemize}
\end{definition}

\begin{remark}
We call it a {\it semicapacity} because given a symplectic embedding  $\psi:
(N,\tau) \to (M,\om)$ such that $\dim N = \dim M$ it cannot be expected, in
general, that $c_{\HZ}^G(N,\tau) \leq c_{\HZ}^{\psi_*G}(M,\om)$, where
$\psi_*: \pi_1(N) \to \pi_1(M)$ is the induced homomorphism on the fundamental
group. However, we can state the following {\it weak monotonicity property}
(for a proof see Section \ref{proofthmA}): given a symplectic embedding
$\psi: (N,\tau)\to (M,\om)$ such that $\dim N = \dim M$ then
$$ c_{\HZ}^{\psi_*^{-1}H}(N,\tau) \leq c_{\HZ}^H(M,\om),$$
for every subgroup $H \subset \pi_1(M)$. In particular, if $\psi_*$ is injective we have that  $c_{\HZ}^G(N,\tau)
\leq c_{\HZ}^{\psi_*G}(M,\om)$.
\end{remark}

Note that obviously,
$$ c_{\HZ}(M,\om) = c_{\HZ}^{\pi_1(M)}(M,\om) \leq c_{\HZ}^G(M,\om), $$
for every subgroup $G \subset \pi_1(M)$. Moreover, it can be show that if
the Hofer-Zehnder $G$-semicapacity is bounded then there are periodic orbits
with homotopy class in $G$ on almost all energy levels for every proper
Hamiltonian \cite{MS}. A symplectic manifold $(M,\om)$ has {\it bounded Hofer-Zehnder $G$-semicapacity}
if for every open subset with compact closure $U \overset{i}{\hra} M$, we have
that
$$ c^{i^{-1}_*G}_{\HZ}(U,\om) < \infty. $$

We remark that the periodic orbits in $M \times D^2(r)$ given  by the
theorems of Hofer-Viterbo and Lu are {\it contractible}. Thus, in terms of
symplectic semicapacities, their theorems state the stronger result that 
$$ c_{\HZ}^0(M \times D^2(r),\om \oplus \om_0) \leq \pi r^2 $$ 
for any geometrically bounded symplectic manifold $(M,\om)$ and $0 < \pi r^2 \leq m(M,\om)$. By
the weak monotonicity property, we can conclude then that
$c_{\HZ}^{\pi_1(S^1)}(M \times T^*S^1,\om\oplus\om_0)$ is bounded, where
$\pi_1(S^1) \subset \pi_1(M \times T^*S^1)$ is the obvious subgroup generated
by the circle in the cylinder $T^*S^1$. In fact, the kernel of the induced
homomorphism on the fundamental group of the symplectic embedding of $M\times
T^*S^1$ into $M \times \R^2$ is given by $\pi_1(S^1)$.

The aim of the first part of this paper is to generalize this result in the
following way: note that $M\times T^*S^1$ has a natural Hamiltonian free
circle action with Hamiltonian given by the angle form $\Phi$ viewed here
as a function on $TS^1 \simeq T^*S^1$. We will consider the following:

\begin{question}
Is the existence of a Hamiltonian free circle action sufficient to ensure the
boundedness of the Hofer-Zehnder semicapacity (with respect to the subgroup of
$\pi_1(M)$ generated by the orbits of the action)?
\end{question}

Before we state our first result we need a previous topological definition
related to free circle actions:

\begin{definition}
Let $\vr: N \times S^1 \to N$ be a free circle action on a connected open subset $N$ of
a manifold $M$. We define the {\it order of the action} $\vr$ as the order of
the cyclic subgroup of $\pi_1(M)$ generated by the homotopy class of the orbits
of $\vr$.
\end{definition}

Now, we are able to state our first theorem. In what follows we will denote
the oscillation of a Hamiltonian $H$ by $\|H\| := \sup H - \inf H$.

\begin{theorem}
\label{thmA}
Let $(M,\om)$ be a geometrically bounded symplectic manifold and $N \subset M$ a connected open subset
that admits a free Hamiltonian circle action $\vr: N \times S^1 \to N$ (with
period equal to $2\pi$) given by the Hamiltonian $H_1: N \to \R$. Suppose that
the order of this action $n_\vr$ (considered as an action on $M$) satisfies
$n_\vr > 2$.\footnote{Recently, we removed this technical hypothesis on the order of the action, see \cite{Mac}.} Given an open subset $U \overset{i}{\hookrightarrow} N$ with
compact closure, suppose that
$$ m(M,\om) > \|H_1|_U\|\bigg(1+\frac{2}{n_\vr-2}\bigg). $$
Then, we have that 
$$ c^{i_*^{-1}G_\vr}_{\HZ}(U,\om) \leq
2\pi\|H_1|_U\|\bigg(1+\frac{2}{n_\vr-2}\bigg) < \infty, $$
where $i_*: \pi_1(U) \to \pi_1(N)$ is the induced homomorphism on the
fundamental group and $G_\vr \subset \pi_1(N)$ is the subgroup generated by the
orbits of $\vr$.
\end{theorem}

\begin{remarks}
\begin{itemize}
\item Note that the previous theorem cannot be stated directly for $N$ because
we need to assume that $M$ is geometrically bounded and it not clear that an open subset of a geometrically bounded symplectic
manifold is geometrically bounded.

\item The circle action above needs to be Hamiltonian. In fact, the example of
Herman-Zehnder \cite{Her,Zeh} has unbounded Hofer-Zehnder capacity and it is
easy to see that it admits a free symplectic circle action with infinite order.
\end{itemize}
\end{remarks}

Before we explain the main ideas involved in the proof of Theorem
\ref{thmA}, let us state some corollaries and applications.

\begin{corollary}
\label{coreulerclass}
Let $(M,\om)$ be a geometrically bounded symplectic manifold. Suppose that $M$ admits a free
Hamiltonian circle action generated by $H_1:M \to \R$ whose Euler class $[e]
\in H^2(M/S^1,\R)$ satisfies $[e]|_{\pi_2(M/S^1)} = 0$. Then, for any open
subset $U \overset{i}{\hra} M$ with compact closure, we have that
$$ c_{\HZ}^{i^{-1}_*G_\vr}(U,\om) \leq 2\pi\|H_1|_U\| < \infty, $$
provided that $m(M,\om) > \|H_1|_U\|$.
\end{corollary}

Actually, it can be proved that given a free circle action on a manifold $M$
such that the Euler class $[e] \in H^2(M/S^1,\R)$ satisfies 
$[e]|_{\pi_2(M/S^1)} = 0$ then the order of the action is infinite (see
Section \ref{eulerclass}). 

Now, consider the finite subgroup $\Z_n \subset S^1$. Note that, since the
circle action is Hamiltonian, the induced action of $\Z_n$ on $M$ is
symplectic. Since the action is free, the quotient manifold $M/\Z_n$ has  a
induced symplectic structure $\om^n$ defined uniquely by the property that
the  pullback of $\om^n$ by the quotient projection is $\om$. Moreover,
$M/\Z_n$ has an induced free Hamiltonian circle action (with period $2\pi/n$)
with order greater than or equal to $n$. Consequently, we have the following
immediate corollary:

\begin{corollary}
\label{cor1}
Let $(M,\om)$ be a geometrically bounded symplectic manifold. Suppose that $M$ admits a free
Hamiltonian circle action generated by $H_1:M \to \R$ and let $U \subset M$ be
an open subset with compact closure. If $n\geq 3$ and $m(M,\om) \geq
(1/n)\|H_1|_U\|$, then
$$ c^{i^{-1}_*G_\vr}_{\HZ}(\tau_n(U),\om) \leq (2\pi/n)\|H_1|_U\|(1+2/(n-2)) < \infty, $$
where $\tau_n: M \to M/\Z_n$ is the quotient projection, $i: \tau_n(U) \to
M/\Z_n$ is the inclusion, $\om$ is the induced symplectic form on $M/\Z_n$ and
$\vr$ is the induced circle action.
\end{corollary}

There are a lot of non-trivial examples of symplectic manifolds admitting such
free Hamiltonian circle actions. The first one that we will consider is the
following: let $M$ be a compact manifold admitting a free circle action. Then
the lift of this action to the cotangent bundle of $M$ is Hamiltonian with
respect to the canonical symplectic form. In fact, it is Hamiltonian with
respect to any twisted symplectic form $\om_0 + \pi^* \Om$, where $\Om$ is a
closed two-form on $M$ given by the pullback of a closed two-form $\Om_{M/S^1}$
on $M/S^1$. Thus, we have the following corollaries:

\begin{corollary}
\label{cor2}
Let $M$ be a compact manifold admitting a free circle action $\vr$ whose Euler
class $[e] \in H^2(M/S^1,\R)$ satisfies $[e]|_{\pi_2(M/S^1)} = 0$ and $\Om$
a weakly exact two-form on $M$ given by the pullback of a two-form
$\Om_{M/S^1}$ on $M/S^1$ . Then the cotangent bundle  $(T^*M,\om_0 +
\pi^*\Om)$ has bounded Hofer-Zehnder $G_\vr$-semicapacity.
\end{corollary}

\begin{corollary}
Let $M$ be a compact manifold such that $\pi_2(M)=0$ and $P$ the total space of
a circle bundle $S^1 \overset{\vr}{\lra} P \lra M$. Then $(T^*P,\om_0)$ has
bounded Hofer-Zehnder $G_\vr$-semicapacity.
\end{corollary}

\begin{corollary}
\label{cor3}
Let $M$ be a compact manifold admitting a free circle action $\vr$. Then
$(T^*(M/\Z_n),\om_0 + \pi^*\Om)$ has bounded Hofer-Zehnder $G_\vr$-semicapacity
for any $n \geq 3$ and any weakly exact 2-form $\Om = \tau_n^*\Om_{M/S^1}$, where
$\tau_n: M/\Z_n \to M/S^1$ is the quotient projection.
\end{corollary}

It is easy to construct specific non-trivial examples of such circle bundles.
For instance, the lens spaces $S^3/\Z_n$ $(n>2)$ over $S^2$, the Heisenberg
manifold over the 2-torus and, more generally, any 2-step nilmanifold
\cite{PS}. In fact, it can be proved that the group of equivalence classes of
circle bundles over a manifold $M$ is isomorphic to $H^2(M,\Z)$ \cite{Kob}.
Thus, any non-trivial element of $H^2(M,\Z)$ for a manifold $M$ such that
$\pi_2(M)=0$ corresponds to a non-trivial circle bundle with infinite order.
The Heisenberg manifold, for example, corresponds to the cohomology class of
the area form on $T^2$.

In the Section \ref{proofthmB}, we show some applications of the Theorem
\ref{thmA} where the circle action is not a lifted action. We will consider
circle actions given by the magnetic flow associated to monopoles on surfaces.

More precisely, let $(M,g)$ be a closed Riemannian manifold, $\Om$ a closed
2-form on $M$ and $TM$ endowed with the twisted symplectic form $\om_0 +
\pi^*\Om$, where $\om_0$ is the pullback of the canonical symplectic form via
the Riemannian metric. When $M$ is a hyperbolic surface and $\Om$ is the area
form, this manifold has remarkable properties and is related to McDuff's
example of a symplectic manifold with disconnected contact-type boundary
\cite{McD} (see Section \ref{proofthmB}, where we give a simple
construction of such a manifold). The Hamiltonian flow $\phi$ given by the
kinetic energy $H(x,v) = (1/2)g_x(v,v)$ with respect to $\om_0 + \pi^*\Om$ is
called the magnetic flow generated by the Riemannian metric $g$ and the
magnetic field $\Om$.

It can be show that, when $M$ is a surface of constant negative curvature and
$\Om$ is the area form, the magnetic flow defines a free circle action with
infinite order on certain subsets of $TM$. Using this fact and the expanding
completion of convex manifolds developed by Eliashberg and Gromov \cite{EG}, we
can prove the following theorem:

\begin{theorem}
\label{thmB}
Let $M$ be a surface of genus $g \geq 2$ endowed with the hyperbolic metric
$g_0$ and the K\"ahler form $\Om_{g_0}$. Let
$H_{g_0}: TM \to \R$ be the usual Hamiltonian given by the kinetic energy
$H_{g_0}(x,v)= (1/2)g_0(v,v)$ and $k \in \R$ a positive real number. Consider
the open subset $U_k = \cup_{0<\mu<k}H_{g_0}^{-1}(\mu)$. Then,
$$ c^G_{\HZ}(U_k,\om_0 + \pi^*\Om_{g_0}) < \infty, $$
as long as $0<k < 1/2$, where $G \subset \pi_1(U_k)$ is the subgroup
generated by the fibers of the unitary bundle. In particular, the periodic orbits have contractible projection on $M$. On the other hand, we have that $c^G_{\HZ}(U_k,\om_0 + \pi^*\Om_{g_0}) = \infty$ for every $k>1/2$.
\end{theorem}

By Moser's theorem, we have the following corollary:
 
\begin{corollary}
Let $M$ be a surface of genus $g\geq 2$ and $\Om$ a  symplectic form on $M$.
Then there exists a Riemannian metric $g$ on $M$ such that
$$ c^G_{\HZ}(U_k,\om_0 + \pi^*\Om) < \infty, $$
as long as $0 < k < 1/2$, where $U_k = \cup_{0<\mu<k}H_g^{-1}(\mu)$ and $G
\subset \pi_1(U_k)$ is the subgroup generated by the fibers of the unitary
bundle as above.
\end{corollary}

\begin{remark}
When $M$ is the 2-torus and $\Om$ is a symplectic form, it is well know
that
$$ c_{\HZ}^{ker\, i_*}(U,\om_0 + \pi^*\Om) < \infty, $$
for every open subset $U \overset{i}{\hra} TT^2$ with compact closure. In fact,
it is easy to prove that $(TT^2,\om_0 + \pi^*\Om)$ is symplectomorphic to
$(T^2\times \R^2,\Om \oplus \sigma)$, where $\sigma$ is the canonical
symplectic form on $\R^2$ \cite{GK}. This result can also be obtained by the
Theorem \ref{thmA}.
\end{remark}

Now, we will use Theorem \ref{thmA} to show that for closed K\"ahler
manifolds $(M,\Om)$ the possible obstruction to the boundedness of the
Hofer-Zehnder capacity is a very thin subset of $M$ which, in many situations,
can be explicitly described.

To prove it, we will need the following nice decomposition result of P. Biran
\cite{Bir1} which enable us to represent a K\"ahler manifold as a disjoint
union of two basic components whose symplectic nature is very standard:

\begin{theorem}[P. Biran \cite{Bir1}]
\label{Bir}
Let $(M^{2n},\Om)$ be a closed K\"ahler manifold with $[\Om] \in
H^2(M,\Z)$ and $\Sigma \subset M$ a complex hypersurface whose homology class
$[\Sigma] \in H_{2n-2}(M)$ is the Poincar\'e dual to $k[\Om]$ for some $k \in
\N$. Then, there exists an isotropic CW-complex $\Delta \subset (M,\Om)$ whose
complement - the open dense subset $(M\sm\Delta,\Om)$ - is symplectomorphic to
a standard symplectic disc bundle $(E_0,\frac{1}{k}\om_{\text{can}})$ modeled on the
normal bundle $N_\Sig$ of $\Sig$ in $M$ and whose fibers have area $1/k$.
\end{theorem}

The symplectic form $\om_{\text{can}}$ is given by
$$ \om_{\text{can}} = k\pi^*(\Om|_\Sigma) + d(r^2\al), $$
where $\pi: E_0 \to \Sigma$ is the bundle projection, $r$ is the radial
coordinate using a Hermitian metric $\|\cdot\|$ and $\al$ is a connection form
on $E$ such that $d\al = -k\pi^*(\Om|_\Sigma)$. The form
$\frac{1}{k}\om_{\text{can}}$ is uniquely characterized by the requirements
that its restriction to the zero section $\Sig$ equals $\Om|_\Sig$, the fibers
of $\pi: E_0 \to \Sigma$ are symplectic and have area $1/k$ and
$\om_{\text{can}}$ is invariant under the obvious circle action along the
fibers. It is called standard because the symplectic type of
$(E_0,\om_{\text{can}})$ depends only on the symplectic type of
$(\Sig,\Om|_\Sig)$ and the topological type of the normal bundle $N_\Sig$
\cite{Bir1}.

Let us recall that the pair $(M,\Sig)$ is called {\it subcritical} \cite{BC1,BC2} if $M\sm\Sig$ is
a subcritical Stein manifold, that is, if there exists a plurisubharmonic Morse
function $\vr$ on $M\sm\Sig$ such that $\text{index}_p(\vr) < \dim_\C M$ for
every critical point $p$ of $\vr$. It is equivalent to the condition that
the dimension of $\Delta$ (that is, the maximal dimension of the cells of
$\Delta$) is strictly less than $n$.

\begin{theorem}
\label{thmC}
Let $(M^{2n},\Om)$ be a closed K\"ahler manifold ($n>2$) with $[\Om] \in
H^2(M,\Z)$ and $\Sigma \subset M$ a complex hypersurface whose homology class
$[\Sigma] \in H_{2n-2}(M)$ is the Poincar\'e dual to $k[\Om]$ for some $k \in
\N$. Then there exists an open dense subset $E$ of $M$ symplectomorphic to a
standard symplectic disc bundle over $\Sigma$ whose complement is an isotropic
CW-complex $\Delta \subset M$ such that if $[\Om]|_{\pi_2(M)}=0$ then
$$ c^G_{\HZ}(E\sm\Sig,\Om) \leq 1/k, $$
where $G \subset \pi_1(E\sm\Sigma)$ is the subgroup generated by the orbits of
the obvious $S^1$-action on $E\sm\Sig$. If $[\Om]|_{\pi_2(M)} \neq 0$
and $k>2$, we have that
$$ c^G_{\HZ}(E\sm\Sig,\Om) \leq \frac{1}{k} + \frac{2}{k^2-2k}. $$
Moreover, if $(M,\Sig)$ is subcritical and $[\Om]|_{\pi_2(M)}=0$, then 
$$ c^{i_*G}_{\HZ}(M\sm\Sigma,\Om) \leq 1/k, $$
where $E\sm\Sig \overset{i}{\hra} M\sm\Sig$.
\end{theorem}

\begin{remarks}
\begin{itemize}
\item In \cite{Mac} it is proved that the hypothesis that $n>2$ and that either $[\Om]|_{\pi_2(M)} = 0$ or $k>2$ are not necessary to ensure that $c^G_{\HZ}(E\sm\Sig,\Om)$ is finite.

\item The result above cannot be stated directly for symplectic disk bundles
and concluded as a property of symplectic disk bundles plus Biran's
decomposition result because symplectic disk bundles are not geometrically bounded. Actually, a
fundamental ingredient in the proof of Theorem \ref{thmC} is that, by the
Lefschetz theorem, the circle action on $E\sm\Sig$ has the same order
considered both as an action on $E\sm\Sig$ itself and as an action on
$M\sm\Sig$ (this construction is not necessary in \cite{Mac} since we do not need there the topological assumption on the order of the action).

\item  Under the hypothesis that either $\dim_\R M \leq 6$ or
$[\Om]|_{\pi_2(M)}=0$, it was proved by P. Biran \cite{Bir1} that the {\it
Gromov capacity} of $E$ satisfies $c_G(E,\Om) \leq 1/k$. These assumptions are
not necessary as was proved recently by G. Lu \cite{Lu3}.

\item The result above for {\it subcritical} manifolds was proved by C. Viterbo
in \cite{Vit2}.

\item When $(M,\Sig)$ is subcritical we do not consider the case $k>2$ because, as
was proved in \cite{EGH}, every subcritical polarization has degree 1.
\end{itemize}
\end{remarks}

We remark that the CW-complex $\Delta$ above is given by the union of the
stable manifolds of the gradient flow of a plurisubharmonic function $\vr$
defined on $M\sm\Sig$ such that it can be explicitly computed in many examples
\cite{Bir1}.

The next theorem shows that we can get a similar result in the neighborhood of
a symplectic hypersurface of Donaldson type \cite{Don}:

\begin{theorem}
\label{thmD}
Let $(M^{2n},\Om)$ be a closed symplectic manifold ($n>2$) with $[\Om] \in
H^2(M,\Z)$ and $\Sigma \subset M$ a symplectic hypersurface whose homology
class $[\Sigma] \in H_{2n-2}(M)$ is the Poincar\'e dual to $k[\Om]$ for some $k
\in \N$. Then there exists $k_0\geq 2$ such that if $k>k_0$ there exists an
open neighborhood $V$ of $\Sigma$ such that if $[\Om]|_{\pi_2(M)}=0$ then
$$ c^G_{\HZ}(V\sm\Sig,\Om) \leq 1/k, $$
where $G \subset \pi_1(V\sm\Sigma)$, as above, is the subgroup generated by the
orbits of the obvious $S^1$-action on $V\sm\Sig$. If
$[\Om]|_{\pi_2(M)} \neq 0$, we have that
$$ c^G_{\HZ}(V\sm\Sig,\Om) \leq \frac{1}{k} + \frac{2}{k^2-2k}. $$
\end{theorem}

It is important to remark that an extension of Biran's result for symplectic
hypersurfaces considered above (which is generally expected to remains valid)
will enable us to extend Theorem \ref{thmC} for general symplectic
manifolds.

The essential idea in the proof of Theorem \ref{thmA} is to relate the
Hofer-Zehnder capacity of a symplectic manifold endowed with a Hamiltonian
circle action with the Hofer-Zehnder capacity of its reduced symplectic
manifold (in the sense of Marsden-Weinstein) with respect to this action.

More precisely, we consider the diagonal Hamiltonian circle action on $P:= M \times
T^*S^1$ whose reduced symplectic manifold $(J^{-1}(\mu)/S^1,\sigma_\mu)$ is
given by $(M,\om)$. Then, given a pre-admissible Hamiltonian $H \in
\H(U,\om)$ we construct a $S^1$-invariant pre-admissible Hamiltonian $\HH$ on
$P$ whose reduced dynamics is given by a reparametrization of the Hamiltonian
vector field of $H$ on $M$. Thus, we can apply the results of Hofer-Viterbo and
Lu to get a periodic orbit for $X_\HH$ and so, by reduction, a periodic  orbit
for $X_H$. The essential step is to show the non-triviality of the projected
closed orbit and it is here that the hypothesis on the order of the action
plays an essential role. The idea is to use the condition on the homotopy of
the periodic orbit together with an upper bound on the period of the orbit to ensure that the orbit cannot be tangent to the trajectories of the diagonal action.

In the second part of this paper, we will consider the problem of the existence
of periodic orbits on prescribed energy levels for a special class of
Hamiltonian dynamical systems given by the magnetic flows. This problem was
first considered by V. Arnold \cite{Ar1,Ar2} and S. Novikov \cite{Nov} and for
a discussion of the results in this area we refer to \cite{CMP,Gin}.

It was proved by L. Polterovich \cite{Pol1} that for every nontrivial weakly
exact magnetic field on a manifold whose Euler characteristic vanishes, there
exist non-trivial contractible closed orbits of the magnetic flow in a sequence
of arbitrarily small energy levels. The proof uses the geometry of the Hofer's
metric in the group of Hamiltonian diffeomorphisms and the fundamental fact
that the displacement energy of the zero section of $TM$, with respect to the
twisted symplectic form given by a non-vanishing magnetic field, is equal to
zero. Recently, E. Kerman \cite{Ker} gives the same result for magnetic fields
given by symplectic forms. In \cite{Mac} we proved the existence of contractible closed orbits for almost all low energy levels, provided that the magnetic field is also symplectic.

We will prove here an extension of Polterovich's theorem for any manifold,
without the assumption on the Euler characteristic.

\begin{theorem}
\label{thmE}
Let $M$ be any closed Riemannian manifold and $\Om$ a non-trivial weakly exact
magnetic field. Then there exists a sequence of arbitrarily small energy levels
containing non-trivial contractible periodic orbits.
\end{theorem}

The main idea in the proof is similar to that of Theorem \ref{thmA}. In
fact, we consider a lift of the magnetic flow to $T^*M \times T^*S^1$ and use
the topological condition on the periodic orbit to ensure the non-triviality of
the projected periodic orbit by symplectic reduction. The details are given in
the Section \ref{proofthmE}.

\medskip
{\it Acknowledgements: } I am very grateful to Gabriel Paternain for many
helpful comments and suggestions and to the Centro de Investigaci\'on en
Matem\'atica, Guanajuato, Mexico for the hospitality while part of this work
was in progress. I am also grateful to Kai Cieliebak, Dusa McDuff and Felix Schlenk for very
useful comments.

\section{Proof of Theorem \ref{thmA}}
\label{proofthmA}

Before we begin the proof of Theorem \ref{thmA} we need some auxiliary
results.

\begin{proposition}
\label{HZin2}
Let $(M,\om)$ be a geometrically bounded symplectic manifold and $P:= M \times T^*S^1$
endowed with the symplectic form $\om_P:= \om \oplus \om^{S^1}_0$, where
$\om^{S^1}_0$ is the canonical symplectic form on $T^*S^1$. Let $U \subset P$
be an open subset with compact closure and suppose that $m(M,\om) \geq
\big|\int_{\tau_2(U)}\om^{S^1}_0\big|$, where $\tau_2: P \to T^*S^1$ is the
projection onto the second factor. Then given a Hamiltonian $H \in
\H(U,\om_P)$, there exists a nonconstant periodic orbit $\gamma$ of $X_H$ with
period
$$ T < T_{max}(U,H) := \frac{\bigg|\int_{\tau_2(U)}\om^{S^1}_0\bigg|}{m(H)}. $$
Moreover,
the homotopy class $[\ga]$ of $\ga$ belongs to the subgroup $\pi_1(S^1)
\subset \pi_1(P)$.
\end{proposition}

\begin{proof}
Since $\ov U$ is compact, there exists a positive constant $a > 0$ such that
$\tau_2(U) \subset S^1 \times [-a/2,a/2] \subset T^*S^1$.  Let $\om_0$ be the
canonical symplectic form on $\R^2$ and consider the symplectomorphism $\phi:
S^1 \times [-a,a] \to A := \{(x,y) \in \R^2; a \leq x^2+y^2 \leq 5a \}$ given
by
$$ \phi(\theta,r) = (\sqrt{3a+2r}\sin\theta,\sqrt{3a+2r}\cos\theta). $$
Note that $Ker(Id,\phi)_*= \pi_1(S^1)$, where $(Id,\phi)_*: \pi_1(T^*M
\times S^1 \times [-a,a]) \simeq \pi_1(T^*P) \to \pi_1(T^*M \times \R^2) \simeq
\pi_1(T^*M)$ is the homomorphism induced on the fundamental group by the
transformation given by the identity and $\phi$ in the first and second factors
respectively.

Now, notice that it is sufficient to prove the Proposition for an open subset $U^\prime$ such that $U \subset U^\prime$ and
$$ \bigg|\int_{\tau_2(U^\prime)}\om^{S^1}_0 - \int_{\tau_2(U)}\om^{S^1}_0\bigg| < \epsilon, $$
for $\epsilon>0$ arbitrarily small.

Thus, we can suppose, without loss of generality, that $\phi(\tau_2(\ov U))
\subset \R^2$ is a connected two-dimensional compact submanifold with boundary,
such that there exists an open disk of radius $R$ with $L$ distinct points ($0
\leq L < \infty$) $y_j \in B^2(R)$ and an orientation preserving diffeomorphism
$$ \psi: \phi(\tau_2(U)) \to B^2(R)\setminus\{y_1,...,y_L\} $$
such that
$$ \bigg|\int_{\phi(\tau_2(U))}\om_0\bigg| = 
\bigg|\int_{B^2(R)\setminus\{y_1,...,y_L\}}\om_0\bigg|
= \pi R^2.$$

From a theorem of Dacorogna and Moser (see Lemma 2.2 in \cite{Si}) this $\psi$
can be required to be symplectic. Thus, consider the Hamiltonian $\bar H: M
\times B^2(R) \to \R$ given by
\begin{equation*}
\bar H(x,y)= 
\begin{cases}
m(H) \text{ if } y=y_j \text{ for some } j = 1,...,L \\
H(x,\phi^{-1}\psi^{-1}y) \text{ otherwise}
\end{cases}
\end{equation*}
which is obviously $C^\infty$, since $H|_{U\setminus K}=m(H)$, where $K \subset
U$ is a compact subset such that $K \subset U\setminus\partial U$.

By Theorem C of \cite{Lu1}, $X_{\bar H}$ has a contractible periodic
orbit $\bar\ga$ with period
$$ T < \frac{\pi R^2}{m(\bar H)} =
\frac{\big|\int_{\phi(\tau_2(U))}\om_0\big|}{m(\bar H)} =
\frac{\big|\int_{\tau_2(U)}\om^{S^1}_0\big|}{m(H)}. $$

Finally, note that, by the remark above, the periodic orbit $\ga$ of $X_H$ given
by $\ga = (Id,\psi\circ\phi)^{-1}\bar\ga$ has homotopy class $[\ga]$ contained
in $\pi_1(S^1) \subset \pi_1(P)$, since $\bar\ga$ is contractible.
\end{proof}

\begin{lemma}
If $\psi: (M_1,\om_1) \to (M_2,\om_2)$ is a symplectic embedding, then
$$ c_{\HZ}^{\psi_*^{-1}(\psi_*G)}(M_1,\om_1) \leq c_{\HZ}^{\psi_*G}(M_2,\om_2), $$
where $G \subset \pi_1(M_1)$ is a subgroup of the fundamental group of $M_1$.
\label{weakmonotonicity}
\end{lemma}

\begin{proof}
Define the map $\psi_*: \H(M_1,\om_1) \to C^\infty(M_2)$ by
\begin{equation*}
\psi_*(H)=
\begin{cases}
H \circ \psi^{-1}(x)\text{ if } x\in\psi(M_1);\\
m(H)\text{ if } x\notin\psi(M_1).\\
\end{cases}
\end{equation*}
It is clear that $\psi_*(\H(M_1,\om_1)) \subset \H(\psi(M_1),\om_2)$ and that
$m(\psi_*H) = m(H)$.

Then it is sufficient to prove that
$\psi_*(\H^{\psi_*^{-1}(\psi_*G)}_a(M_1,\om_1)) \subset
\H^{\psi_*G}_a(\psi(M_1),\om_2)$. To do it, suppose, by contradiction, that there
exists $H \in \H^{\psi_*^{-1}(\psi_*G)}_a(M_1,\om_1)$ such that
$X^{\om_2}_{\psi_*H}$ has a nonconstant periodic orbit $\gamma$ such that
$[\gamma] \in \psi_*G$ and the period of $\gamma$ is less than 1. Then the
image of $\gamma$ by $\psi^{-1}$ (note that the image of $\gamma$ is contained
in $\psi(M_1)$, since $\psi_*H \equiv m(H)$ outside $\psi(M_1)$) is a periodic
orbit of $X^{\om_1}_H$ of the same period and 
$$ [\psi^{-1} \circ \gamma] \in \psi_*^{-1}(\psi_*G), $$
a contradiction.
\end{proof}

\begin{proof}[Proof of Theorem \ref{thmA}: ]
We will divide the proof in three steps:

\subsection{Construction of a reparametrized diagonal circle action
$\rho$ on $N \times T^*S^1$ and its symplectic reduction.}

Let $P = N \times T^*S^1 \subset M \times T^*S^1$, where $M \times T^*S^1$ is
endowed with the symplectic form $\om_P = \om \oplus \om^{S^1}_0$, where
$\om^{S^1}_0$ is the canonical symplectic form on $T^*S^1$. Note that $T^*S^1$
has a natural Hamiltonian circle action given by the Hamiltonian
$\Phi(\theta,\mu)=\mu$ . Let $n \in \R$ be such that $n>2\|H_1|_U\|/(n_\vr-2)$
and $m(M,\om) \geq \|H_1|_U\|+n$.

Consider the diagonal action $\rho: P \times S^1 \to P$ on $P$ given by the
Hamiltonian
\begin{align*}
J(x,\theta,\mu) & = (1/n)(H_1(x) + \mu) \\
& = (1/n)(H_1(x) + \Phi(\theta,\mu)),
\end{align*}
and let $\psi: P \to P$ be the diffeomorphism given by
$$ \psi(x,\theta,\mu) = (\vr_\theta(x),\theta,\mu - H_1(x)), $$
where $\vr: N \times S^1 \to N$ is the flow generated by $H_1$. The factor $n$ is not important in this section and it will be used only in the next section in the construction of a $\rho$-invariant pre-admissible Hamiltonian from $H$.

The following lemma is a straightforward computation.

\begin{lemma}
\label{symptriv}
The diffeomorphism $\psi$ is a symplectomorphism with respect to $\om_P$.
Moreover, $J \circ \psi = (1/n)\Phi$.
\end{lemma}

\begin{proof}
It is clear that $\psi^*J = (1/n)\Phi$. In effect,
\begin{align*}
J(\psi(x,\te,\mu)) & = (1/n)(H_1(\vr_\te(x)) + \mu - H_1(x)) \\
& = (1/n)(H_1(x) - H_1(x) + \mu) \\
& = (1/n)\mu.
\end{align*}

To show that $\psi^*\om_P = \om_P$, it is more convenient to write $\psi$ as 
$$ \psi(z) = \vr_{\theta(z)}(z) - H_1(z)Y, $$
where $\theta: P \to S^1$ is the projection onto the circle and
$Y$ is the unit vector field tangent to the fibers of $T^*S^1$ such that
$\Phi(Y)=1$. Consequently,
$$ d\psi(z)\xi = (d\vr)_{\te(z)}(z)\xi + \pi_2^*d\te(\xi)X_{H_1}(\vr_{\te(z)}(z)) -
dH_1(z)\xi Y(\vr_{\te(z)}(z)), $$
where $X_{H_1}$ is the Hamiltonian vector field generated by $H_1$ and
$\pi_2^*d\te$ is the pullback of the angle form $d\te$ on $S^1$ to $P$. Thus,
\begin{align*}
& (\psi^*\om_P)_z(\xi,\eta) = (\om_P)_{\psi(z)}(d\psi(z)\xi,d\psi(z)\eta) \\
& = \om_P((d\vr)_{\te(z)}(z)\xi,(d\vr)_{\te(z)}(z)\eta) 
+ \om_P((d\vr)_{\te(z)}(z)\xi,\pi_2^*d\te(\eta)X_{H_1}(\vr_{\te(z)}(z))) \\
& - \om_P((d\vr)_{\te(z)}(z)\xi,dH_1(z)\eta Y(\vr_{\te(z)}(z)))
+ \om_P(\pi_2^*d\te(\xi)X_{H_1}(\vr_{\te(z)}(z)),(d\vr)_{\te(z)}(z)\eta) \\
& + \om_P(\pi_2^*d\te(\xi)X_{H_1}(\vr_{\te(z)}(z)),\pi_2^*d\te(\eta)X_{H_1}(\vr_{\te(z)}(z)))\\
& -\om_P(\pi_2^*d\te(\xi)X_{H_1}(\vr_{\te(z)}(z)),dH_1(z)\eta Y(\vr_{\te(z)}(z)))\\
& - \om_P(dH_1(z)\xi Y(\vr_{\te(z)}(z)),(d\vr)_{\te(z)}(z)\eta)
-\om_P(dH_1(z)\xi Y(\vr_{\te(z)}(z)),\pi_2^*d\te(\eta)X_{H_1}(\vr_{\te(z)}(z)))\\
& + \om_P(dH_1(z)\xi Y(\vr_{\te(z)}(z)),dH_1(z)\eta Y(\vr_{\te(z)}(z))).
\end{align*}
But note that
$$\om_P(\pi_2^*d\te(\xi)X_{H_1}(\vr_{\te(z)}(z)),\pi_2^*d\te(\eta)X_{H_1}(\vr_{\te(z)}(z)))=0$$
and
$$\om_P(dH_1(z)\xi Y(\vr_{\te(z)}(z)),dH_1(z)\eta Y(\vr_{\te(z)}(z)))=0,$$
since the vectors are colinear. On the other hand,
$$\om_P(\pi_2^*d\te(\xi)X_{H_1}(\vr_{\te(z)}(z)),dH_1(z)\eta Y(\vr_{\te(z)}(z)))=0$$
and
$$\om_P(dH_1(z)\xi Y(\vr_{\te(z)}(z)),\pi_2^*d\te(\eta)X_{H_1}(\vr_{\te(z)}(z)))=0,$$
because the vectors are orthogonal with respect to the product decomposition $P
= M \times T^*S^1$ and by the definition of $\om_P$.

We have then that,
\begin{align*}
& (\psi^*\om_P)_z(\xi,\eta)
= \om_P((d\vr)_{\te(z)}(z)\xi,(d\vr)_{\te(z)}(z)\eta) 
+ \om_P((d\vr)_{\te(z)}(z)\xi,\pi_2^*d\te(\eta)X_{H_1}(\vr_{\te(z)}(z)))\\
&- \om_P((d\vr)_{\te(z)}(z)\xi,dH_1(z)\eta Y(\vr_{\te(z)}(z)))
+ \om_P(\pi_2^*d\te(\xi)X_{H_1}(\vr_{\te(z)}(z)),(d\vr)_{\te(z)}(z)\eta)\\
&-\om_P(dH_1(z)\xi Y(\vr_{\te(z)}(z)),(d\vr)_{\te(z)}(z)\eta)=\\
& = (\om_P)_z(\xi,\eta) +
\pi_2^*d\te(\eta)(i_{X_{H_1}(\vr_{\te(z)}(z))}\om)(d\vr_{\te(z)}(z)\xi)
- dH_1(z)\eta(i_{Y(\vr_{\te(z)}(z))}\om)(d\vr_{\te(z)}(z)\xi)\\
& + \pi_2^*d\te(\xi)(i_{X_{H_1}(\vr_{\te(z)}(z))}\om)(d\vr_{\te(z)}(z)\eta)
 - dH_1(z)\xi(i_{Y(\vr_{\te(z)}(z))}\om)(d\vr_{\te(z)}(z)\eta),
\end{align*}
where the last equality follows from the fact that $\vr_t^*\om_P=\om_P$. Now,
note that
$$ \pi_2^*d\te(\eta)(i_{X_{H_1}(\vr_{\te(z)}(z))}\om)(d\vr_{\te(z)}(z)\xi)
= \pi_2^*d\te(\eta)dH_1(z)\xi, $$
because $i_{X_{H_1}}\om_P = dH_1$ and $\vr^*dH_1 = dH_1$. On the other hand, we
have that
$$ dH_1(z)\eta(i_{Y(\vr_{\te(z)}(z))}\om)(d\vr_{\te(z)}(z)\xi)
= dH_1(z)\eta \pi_2^*d\te(\xi), $$
because $i_{Y}\om_P = \pi_2^*d\te$ and $\vr^*\pi_2^*d\te = \pi_2^*d\te$ (note
that $\pi_2^*d\te$ here is the pullback of the angle form $d\te$ on $S^1$ to $P$).
Consequently,
\begin{align*}
(\psi^*\om_P)_z(\xi,\eta) & = (\om_P)z(\xi,\eta) + \\
& + \pi_2^*d\te(\eta)dH_1(z)\xi
- dH_1(z)\eta \pi_2^*d\te(\xi)
+ \pi_2^*d\te(\xi)dH_1(z)\eta
- dH_1(z)\xi \pi_2^*d\te(\eta) \\
& = \om_z(\xi,\eta),
\end{align*}
as desired.
\end{proof}

\begin{lemma}
\label{MWreduction}
The Marsden-Weinstein reduced symplectic space at $J^{-1}(\mu)$ is given by $N$
with the reduced symplectic form $\sigma_\mu$ equal to $\om$ for every $\mu
\in \R$.
\end{lemma}

\begin{proof}
We claim that the quotient projection $\pi_\mu: J^{-1}(\mu) \to N$ is given by $\pi_1 \circ \psi^{-1}|_{J^{-1}(\mu)}$, where
$\pi_1: P \to N$ is the projection onto the first factor. In effect, note that,
by the previous lemma, $\psi$ sends the orbits of $X_J$ to the orbits of
$X_{(1/n)\Phi}$. On the other hand, $\pi_1$ is the quotient projection at
$((1/n)\Phi)^{-1}(\mu) = \psi^{-1}(J^{-1}(\mu))$ with respect to the trivial
circle bundle given by the orbits of $X_{\Phi}$.

To show that $\sigma_\mu$ is equal to $\om$, note that
\begin{align*}
\pi_\mu^*\om & = (\psi^{-1})^*i_{((1/n)\Phi)^{-1}(\mu)}^*\pi_1^*\om \\
& = (\psi^{-1})^*i_{((1/n)\Phi)^{-1}(\mu)}^*\om_P \\
& = i_{\psi(((1/n)\Phi)^{-1}(\mu))}^*\om_P \\
& = i_{J^{-1}(\mu)}\om_P,
\end{align*}
where the third equality follows from the fact that $\psi^*\om_P=\om_P$.
\end{proof}

\subsection{Construction of a $\rho$-invariant pre-admissible
Hamiltonian $\HH$ on $N \times T^*S^1$ from $H$.}

In view of the last lemma, consider the map $\ppi: P \to N$ defined by
$$ \ppi = \pi_1 \circ \psi^{-1}, $$
such that $\ppi|_{J^{-1}(\mu)} = \pi_\mu$ is the quotient projection.
Now, let $H$ be a pre-admissible Hamiltonian on $U$.  We will construct from
$H$ a $\rho$-invariant Hamiltonian defined on $P$. Firstly, fix a sufficiently small
constant $\delta>0$ satisfying
$$ n > (2+2\delta)\frac{\|H_1|_U\|}{n_\vr-(2+2\delta)}. $$
Define the new Hamiltonian by
$$ \HH(z) = (H + \al(J(z))(m(H)-H))(\ppi(z)), $$
where $\al: \R \to \R$ is a $C^\infty$ function such that $0 \leq \al \leq 1$,
$\al(\mu) = 1\ \ \forall \mu \notin (\delta,1-\delta)$, $\al(\mu)= 0\ \ \forall
\mu \in (\frac{1}{2}-\delta,\frac{1}{2}+\delta)$ and $|\al^\pr(\mu)| \leq
2+2\delta$ for every $\mu \in [0,1]$ (see figure \ref{graph1}). Thus, $\HH|_{J^{-1}(\mu)}$ is the lift of 
$H_\mu:= H + \al(\mu)(m(H)-H) = (1-\al(\mu))H + \al(\mu)m(H)$ by the quotient
projection $\pi_\mu$.

\begin{figure}[ht]
\label{fctalpha}
\begin{center}
\includegraphics[width=2in]{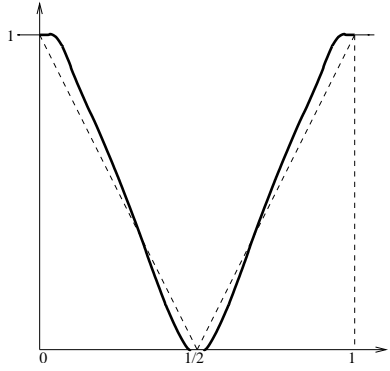}
\caption{\label{graph1} Graph of the function $\al$.}
\end{center}
\end{figure}

The first obvious property of $\HH$ is that $m(\HH) = m(H)$. In effect, $m(\HH)
= \sup_{\mu \in \R} m(H+\al(\mu)(m(H)-H)) = m(H)$, since $0 \leq \al \leq 1$
and $\al(\mu) = 1\ \ \forall \mu \notin (\delta,1-\delta)$.

It is easy to see that $\HH \in \H(\U,\om)$ where $\U$ is given by (see figure \ref{graph2})
$$ \U = \bigcup_{0\leq\mu\leq 1} \pi_\mu^{-1}(U).$$
In effect, it is clear that $0 \leq \HH \leq m(\HH)$. Moreover, $\HH|_{\V}
\equiv 0$ where $\V$ is the subset of $P$ given by
$\bigcup_{\frac{1}{2}-\delta\leq\mu\leq\frac{1}{2}+\delta} \pi_\mu^{-1}(V)$ and
$V \subset M$ is an open set such that $H|_V \equiv 0$ (that exists because $H
\in \H(U,\om)$) since $\al(\mu) = 0\ \ \forall \mu \in
(\frac{1}{2}-\delta,\frac{1}{2}+\delta)$. Note that the interior of $\V$ is not
empty.

On the other hand, $\HH \equiv m(\HH)$ outside the compact set $\K =
\bigcup_{\delta\leq\mu\leq 1-\delta} \pi_\mu^{-1}(K) \subset
\U\setminus\partial\U$, where $K \subset U\setminus\partial U$ is the compact
set such that $H|_{U\setminus K} \equiv m(H)$.

\subsection{Existence of a non-trivial periodic orbit of $X_H$ given
by the reduction of a periodic orbit of $X_\HH$.}

Let us now prove that there exists a periodic orbit $\ga$ of $X_H$ of period
less than $2\pi(\|H_1|_U\|+n)/m(H)$ and homotopy class $[\ga] \in G_\vr$.
Note firstly that 
$$ \bigg|\int_{\tau_2(\U)}\om_0^{S^1}\bigg| = 2\pi (\|H_1|_U\| + n). $$
In effect, the image of $\U \cup J^{-1}(0)$ under the projection $N \times S^1 \times \R \to N \times \R$ is the graph of $-H_1$ restricted to $U$
and the projection of $\U \cup J^{-1}(1)$ is the graph of $n-H_1|_U$:

\begin{figure}[h]
\begin{center}
\psfrag{mu}{$\mu$}
\psfrag{UU}{$\U$}
\psfrag{grH1U}{$-H_1|_U$}
\psfrag{Phi1}{$\Phi^{-1}(n)$}
\psfrag{Phi0}{$\Phi^{-1}(0)$}
\psfrag{J1}{$J^{-1}(1)$}
\psfrag{J0}{$J^{-1}(0)$}
\includegraphics[width=2in]{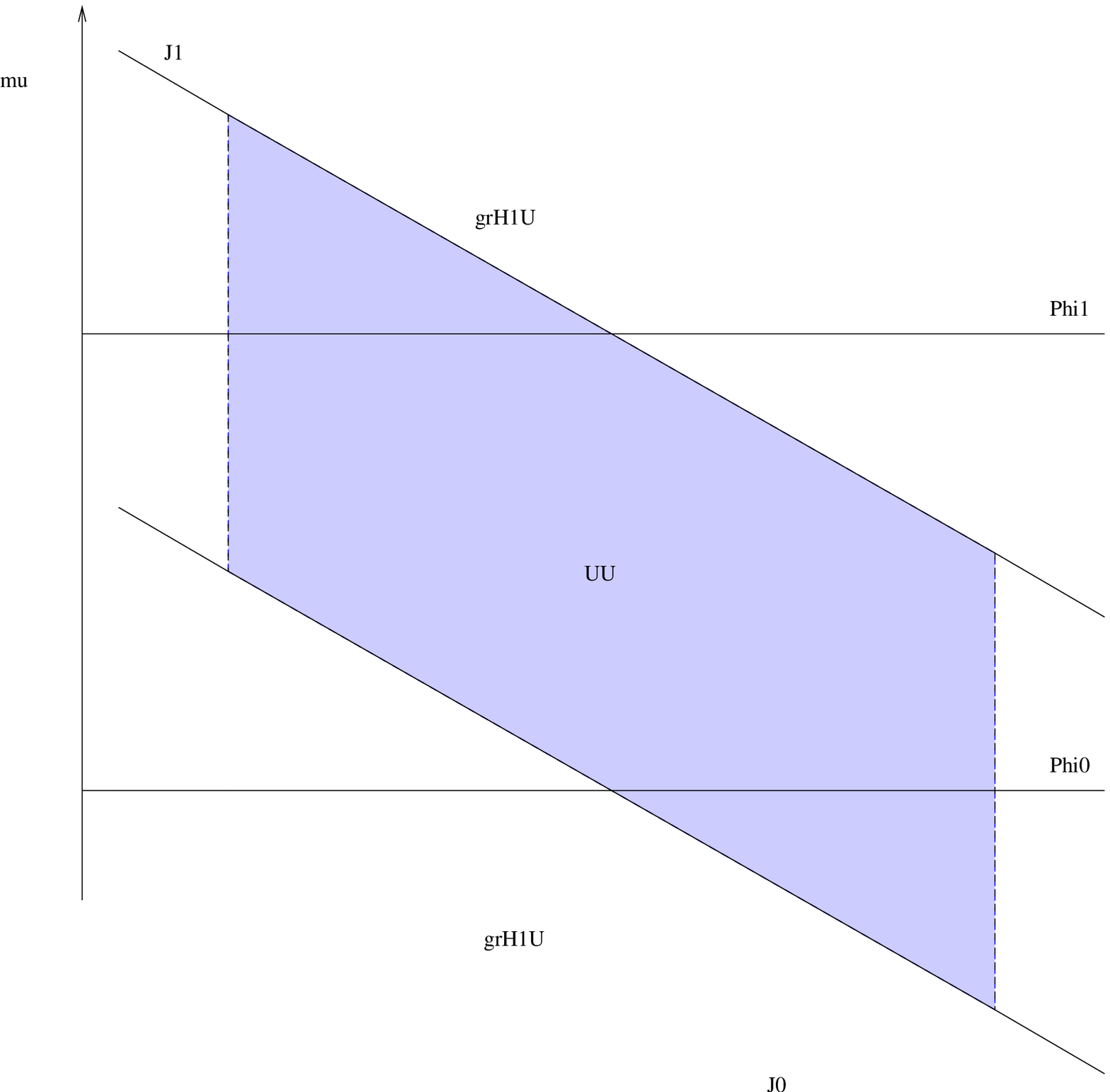}
\caption{\label{graph2} The subset $\U$.}
\end{center}
\end{figure}

\vskip .3cm
\noindent {\bf Reduced dynamics of $X_\HH$:}
\vskip .2cm

Note that $H_\mu = (1-\al(\mu))H + \al(\mu)m(H)$ has the Hamiltonian vector
field with respect to $\om$ given by $X_{H_\mu} = (1-\al(\mu))X_H$. Thus, the
reduced dynamics of $X_\HH$ at $J^{-1}(\mu)$ is a reparametrization of the
dynamics of $X_H$. Since $0 \leq \al(\mu) \leq 1$, the nonconstant periodic
orbits of $X_{H_\mu}$ have period greater than or equal to those of $X_H$.

Consequently, by the Lemma \ref{MWreduction}, it is sufficient to show
that the periodic orbit of $\HH$ given by the Proposition \ref{HZin2} is not
tangent to a fiber given by the action of $\vr$ (such that it is projected
onto a singularity of $X_H$).

\vskip .3cm
\noindent {\bf Non-triviality of the projected orbit:}
\vskip .2cm

Suppose that there exists a periodic orbit $\ga: [0,T] \to \U$ of $X_\HH$ tangent to an orbit
of $\rho$ and whose homotopy class belongs to $\pi_1(S^1) \subset \pi_1(P)$. We
will show that the period $T$ of $\ga$ is strictly greater than $T_{\max}(\U,\HH)$
(see Proposition \ref{HZin2}).

Firstly, note that
$$ T_{\max}(\U,\HH) = \frac{\big|\int_{\tau_2(\U)}\om_0^{S^1}\big|}{m(\HH)} =
\frac{2\pi(\|H_1|_U\|+n)}{m(H)}. $$

Now, let $W$ be an almost complex structure on $P$ compatible with $\om_P$ and
$\lg\cdot,\cdot \rg = \om_P(W\cdot,\cdot)$ be the corresponding Riemannian
metric. We have that along $\ga$,
$$ X_\HH = \bigg\lg X_\HH,\frac{X_J}{\|X_J\|^2}\bigg\rg X_J. $$
On the other hand, the period of the orbits of $X_J$ is equal to $2\pi n$ because the action on $M$ is free.
Consequently, since $[\ga] \in \pi_1(S^1)$, the projection $\pi_1\circ\gamma$ must wrap around itself at least $n_\vr$ times. Hence, $\gamma$ also must wrap around itself at least $n_\vr$ times and so the period of $\ga$ is given by
$$ T = \frac{2\pi n n_\vr}{\big|\big\lg X_\HH,\frac{X_J}{\|X_J\|^2}\big\rg\big|}. $$
But, by the definition of $n$,
$$ T = \frac{2\pi n n_\vr}{\big|\big\lg X_\HH,\frac{X_J}{\|X_J\|^2}\big\rg\big|} >
\frac{2\pi (2+2\delta)(\|H_1|_U\|+n)}{\big|\big\lg X_\HH,\frac{X_J}{\|X_J\|^2}\big\rg\big|}. $$
Actually, note that the inequality
$$ 2\pi nn_\vr > 2\pi (2+2\delta)(\|H_1|_U\|+n) $$
is equivalent to
$$ 2\pi n(n_\vr-(2+2\delta)) > 2\pi(2+2\delta)\|H_1|_U\|, $$
that is,
$$ n > (2+2\delta)\frac{\|H_1|_U\|}{n_\vr-(2+2\delta)}, $$
which is exactly the condition on $n$ and $\delta$. Thus, it is sufficient to prove that
$$ \bigg|\bigg\lg X_\HH,\frac{X_J}{\|X_J\|^2}\bigg\rg\bigg| \leq
(2+2\delta)m(H), $$
which implies
$$ T > \frac{2\pi (\|H_1|_U\| + n)}{m(H)}. $$
To conclude this inequality, note that
\begin{align*}
\bigg|\bigg\lg X_\HH,\frac{X_J}{\|X_J\|^2} \bigg\rg\bigg|
& = \bigg|\bigg\lg W\nabla\HH,\frac{W\nabla J}{\|W\nabla J\|^2} \bigg\rg\bigg| \\
& = \bigg|\bigg\lg \nabla\HH,\frac{\nabla J}{\|\nabla J\|^2} \bigg\rg\bigg| \\
& = \bigg|d\HH\bigg(\frac{\nabla J}{\|\nabla J\|^2}\bigg)\bigg|
\end{align*}
where in the second equation we used the fact that $W$ defines an
isometry. On the other hand, we have that for every $\xi \in T_zP$,
\begin{align*}
d\HH(z)\xi & =
dH(\ppi(z))d\ppi(z)\xi
+ m(H)\al^\pr(J(z))dJ(z)\xi
 - H(\ppi(z))\al^\pr(J(z))dJ(z)\xi\\
& \ \ - \al(J(z))dH(\ppi(z))d\ppi(z)\xi\\
& = (m(H) - H\circ\ppi(z))\al^\pr(J(z))dJ(z)\xi
+ (1-\al(J(z)))dH(\ppi(z))d\ppi(z)\xi.
\end{align*}
But note that, since $\ga$ is tangent to a fiber, it is projected onto a
singularity of $X_H$, that is, $dH(\ppi(z)) = 0$ for every $z \in \ga$.
Consequently, we have that
\begin{align*}
\bigg|\bigg\lg X_\HH,\frac{X_J}{\|X_J\|^2} \bigg\rg\bigg|
& = \bigg|d\HH\bigg(\frac{\nabla J}{\|\nabla J\|^2}\bigg)\bigg| \\
& = (m(H) - H\circ\ppi(x))|\al^\pr(J(z))| \\
& \leq (2+2\delta)m(H),
\end{align*}
as desired.

Finally, note that the periodic orbit $\bar\ga$ of $X_\HH$ given by the
Proposition \ref{HZin2} satisfies $[\bar\ga] \in \pi_1(S^1)$. Consequently,
the corresponding reduced periodic orbit $\ga$ satisfies $[\ga] \in G_\vr$ since
$(\pi_\mu)_*\pi_1(S^1) = G_\vr$.

Thus, we obtain the estimate
$$ c_{\HZ}^{i^{-1}_*G_\vr}(U,\om) \leq 2\pi(\|H_1|_U\|+n). $$
But, since $n$ is arbitrarily close to $2\|H_1|_U\|/(n_\vr-2)$, we conclude that
$$ c_{\HZ}^{i^{-1}_*G_\vr}(U,\om) \leq 2\pi\|H_1|_U\|\bigg(1+\frac{2}{n_\vr-2}\bigg). $$
\end{proof}

\section{Proof of Corollary \ref{coreulerclass}}
\label{eulerclass}
The proof of Corollary \ref{coreulerclass} follows immediately from the
Theorem \ref{thmA} and the following useful proposition:

\begin{proposition}
\label{propeulerclass}
Let $S^1\overset\vr\lra P \overset{\pi_P}\lra M$ be a principal circle bundle
whose Euler class $[\Om] \in H^2(M,\R)$ satisfies $[\Om]|_{\pi_2(M)}=0$. Then
$\vr$ has infinite order.
\end{proposition}

\begin{proof}
Initially, note that $[\Om]|_{\pi_2(M)}=0$ iff $\tau_M^*[\Om] = 0$, where
$\tau_M: \wt{M} \to M$ is the universal cover. In effect, it is clear that if
$\tau_M^*\Om$ is exact, then $[\Om]|_{\pi_2(M)}=0$. To prove the converse, note
that, by the Hurewicz's theorem, we have that $H_2(\wt{M},\Z) \simeq
\pi_2(\wt{M}) \simeq \pi_2(M)$. Thus, the condition that $\Om$ is weakly exact
is equivalent that the integral of $\tau_M^*[\Om]$ over any 2-cycle is equal to
zero. But, by de Rham's theorem, it is equivalent to $\tau_M^*\Om$ be exact.

Now, suppose that the action $\vr$ has finite order. Consider the lift
$\wt{\vr}$ of $\vr$ to the universal cover $\wt{P}$ of $P$.

\medskip
\noindent {\bf Claim 1.} $\wt{P}/\wt{\vr} = \wt{M}$, where $\wt{P}/\wt{\vr}$ is
the quotient of $\wt{P}$ by the action $\wt{\vr}$ and $\wt{M}$ is the universal
cover of $M$.
\medskip

To prove the claim, note that, since the fiber is path connected,
$\wt{P}/\wt{\vr}$ is simply connected. Thus, it is sufficient to prove that
$\wt{P}/\wt{\vr}$ covers $M$. Now, observe that
$$ M = \wt{P}/\Gamma_P/\vr, $$
where $\Gamma_P$ is the group of deck transformations on $\wt{P}$. On the other hand,
$\wt{\vr}$ commutes with the action of $\Gamma_P$, since it is a lifted action
from $P$. Consequently, it induces an action of a discrete group $\Gamma_M$ on
$\wt{P}/\wt{\vr}$ such that
$$ M = \wt{P}/\Gamma_P/\vr = \wt{P}/\wt{\vr}/\Gamma_M, $$ 
proving the claim.

\medskip
\noindent {\bf Claim 2.} The following diagram is commutative:
\begin{equation*}
\begin{CD}
\wt{P} @>\tau_P>> P \\
@V\pi_{\wt{P}} VV @VV\pi_P V \\
\wt{M} @>\tau_M>> M
\end{CD}
\end{equation*}
\medskip

In fact, note that
$$ \pi_P \circ \tau_P: \wt{P} \to \wt{P}/\Gamma_P/\vr $$
and
$$ \tau_M \circ \pi_{\wt{P}}: \wt{P} \to \wt{P}/\wt{\vr}/\Gamma_M. $$
But, as remarked above, $\wt{\vr} \circ \ga_P = \ga_P \circ \wt{\vr}\
\forall\ga_P \in \Gamma_P$ and $\Gamma_M$ is the induced action on
$\wt{P}/\wt{\vr}=\wt{M}$ such that $\pi_P \circ \tau_P = \tau_M \circ
\pi_{\wt{P}}$.

\medskip
\noindent {\bf Claim 3.} the Euler class of $S^1\overset{\wt{\vr}}\lra \wt{P}
\overset{\pi_{\wt{P}}}\lra \wt{M}$ is given by $[\tau_M^*\Om]$, where $\Om$ is the
curvature form of $S^1\overset\vr\lra P \overset{\pi_P}\lra M$.
\medskip

In effect, let $\al$ be a connection form of $\vr$. Then $\tau_P^*\al$ is a
connection for $\wt{\vr}$ and
\begin{align*}
d\tau_P^*\al & = \tau_P^*d\al \\
& = \tau_P^*\pi_P^*\Om \\
& = \pi_{\wt{P}}^*\tau_M^*\Om.
\end{align*}

Thus, if $\tau_M^*\Om$ is exact, we have that $S^1\overset{\wt{\vr}}\lra \wt{P}
\overset{\pi_{\wt{P}}}\lra \wt{M}$ defines a trivial circle bundle, contradicting the fact
that $\wt{P}$ is simply connected.
\end{proof}

\section{Proof of Theorem \ref{thmB}}
\label{proofthmB}
Let $M$ be a surface of genus $g\geq 2$ endowed with a Riemannian metric $g$ of
constant negative curvature. Let $TM$ be the tangent bundle of $M$ with the
symplectic form $\om_0$ given by the pullback of the canonical symplectic form
of $T^*M$ by the bundle isomorphism induced by $g$.

Consider on $M$ the almost complex structure $J$ such that $\Om =
g(J\cdot,\cdot )$ is the K\"ahler form. As defined in the introduction, the magnetic flow associated to $g$ and $\Om$ is
the Hamiltonian flow given by the kinetic energy $H(x,v) = (1/2)g_x(v,v)$ with
respect to the twisted symplectic form $\om_1:= \om_0 + \pi^*\Om$. It is easy
to see that a curve $t \mapsto (\ga(t),\dot{\ga}(t))$ defines an orbit of this
flow if and only if $\nabla_{\dot{\ga}(t)}\dot{\ga}(t) = J(\dot{\ga}(t))$. In
other words, the trajectories are the curves with constant geodesic curvature.

Consequently, the dynamics of this magnetic flow can be explicitly described:
on the energy levels greater than one half, the magnetic flow is topologically
equivalent to the geodesic flow (which has periodic orbits in every nontrivial
homotopy class and does not have projected contractible periodic orbits,
because there is no conjugate points); on the sphere bundle, it coincides with
the horocycle flow and hence is minimal, that is, all the orbits are dense; and
on the energy levels below than one half, all the orbits are closed and homotopic to the fibers of the sphere bundle.

In particular, it defines (after a suitable reparametrization) a Hamiltonian circle action with infinite order on
$U_{1/2}$ (since $\pi_2(M)=0$) such that we can apply Theorem \ref{thmA}.
The problem is the condition on $TM\setminus M_0$ to be geometrically bounded. To circumvent
this problem we need the following proposition:

\begin{proposition}
\label{tame}
Given $\ep > 0$ there exists an exact geometrically bounded symplectic form $\rho_\ep$ on
$TM\setminus M_0$ such that $\rho_\ep$ restricted to $(TM\setminus M_0)
\setminus \ov U_\ep$ coincides with $\om_1$.
\end{proposition}

\begin{proof}
The proof relies on the method of expanding completion of convex manifolds
introduced by Eliashberg and Gromov \cite{EG}. Firstly, we need the following
key lemma which has intrinsic interest:

\begin{lemma}
The vector field $X$ on $TM\setminus M_0$ given by (with respect to the
horizontal and vertical subbundles)
$$ X(x,v)=\bigg(\frac{1}{\|v\|^2}Jv,\bigg(1-\frac{1}{\|v\|^2}\bigg)v\bigg) $$
is conformally symplectic with respect to $\om_1$, that is, $L_X \om_1 = \om_1$.
\end{lemma}

\begin{proof}
To fix the notation, let $\om_0(\xi,\eta) = \lg \pi_*\xi,K\eta \rg - \lg K\xi,
\pi_*\eta \rg$, where $K$ is the curvature map, be the pullback of the
canonical symplectic form of $T^*M$ by the Riemannian metric and
$\pi^*\Om(\xi,\eta) = \lg J\pi_*\xi,\pi_*\eta \rg$. Define $X_1(x,v) =
((1/\|v\|^2)Jv,0)$ and $X_2(x,v) = (0,(1-1/\|v\|^2)v)$. Consider the Liouville
form
$$ \al_{(x,v)}(\xi) = \lg v,\pi_*\xi \rg $$
and the 1-form
$$ \beta_{(x,v)}(\xi) = \lg Jv,K\xi \rg. $$
Note that $d\al = -\om_0$ and $d((1/\|v\|^2)\beta) = \pi^*\Om$, where the last
equality follows from the fact that $(1/\|v\|^2)\beta$ defines a connection form
on $H^{-1}((1/2)\|v\|^2)$ and that $M$ is a hyperbolic surface.

Let us compute the Lie derivative of $\om_1 = \om_0 + \pi^*\Om$ with respect to
$X_1$. Initially, we have that,
\begin{align*}
i_{(Jv,0)} \om_1(\xi) & = \lg Jv,K\xi \rg - \lg v,\pi_*\xi \rg \\
& = \beta - \al
\end{align*}
such that
\begin{align*}
L_{(Jv,0)}\om_1 & = d(i_{(Jv,0)} \om_1)\\
& = d(\beta - \al).
\end{align*}
But, we have that
\begin{align*}
\pi^*\Om & = d((1/\|v\|^2)\beta) \\
& = \frac{1}{\|v\|^2}d\beta - \frac{2}{\|v\|^3} dv \wedge \beta.
\end{align*}
Consequently,
$$ d\beta = \|v\|^2\pi^*\Om + \frac{2}{\|v\|}dv \wedge \beta. $$
Thus,
$$ d(i_{(Jv,0)} \om_1) = \om_0 + \|v\|^2\pi^*\Om + \frac{2}{\|v\|}dv \wedge
\beta $$
and so
\begin{align*}
d(i_{(1/\|v\|^2)(Jv,0)} \om_1) & = \frac{1}{\|v\|^2}\om_0 + \pi^*\Om 
+ \frac{2}{\|v\|^3} dv \wedge \beta
- \frac{2}{\|v\|^3} dv \wedge i_{(Jv,0)} \om_1 \\
& = \frac{1}{\|v\|^2}\om_0 + \pi^*\Om 
+ \frac{2}{\|v\|^3} dv \wedge \beta
- \frac{2}{\|v\|^3} dv \wedge (\beta-\al) \\
& = \frac{1}{\|v\|^2}\om_0 + \pi^*\Om + \frac{2}{\|v\|^3} dv \wedge \al.
\end{align*}
On the other hand, the Lie derivative of $\om_1$ with respect to $X_2$ is given
by
\begin{align*}
d(i_{(0,(1-(1/\|v\|^2))v)}\om_1) & = d(-\al + \frac{1}{\|v\|^2}\al) \\
& = \om_0 - \frac{1}{\|v\|^2}\om_0 - \frac{2}{\|v\|^3}dv \wedge \al
\end{align*}
since
$$ i_{(0,(1-(1/\|v\|^2))v)}\om_1 = -\al + \frac{1}{\|v\|^2}\al. $$
Consequently, we have that
\begin{align*}
L_X \om_1 & = d(i_X\om_1) \\
& = d(i_{X_1+X_2}\om_1) \\
& = \frac{1}{\|v\|^2}\om_0 + \pi^*\Om +
\frac{2}{\|v\|^3}dv \wedge \al + \om_0 - \frac{1}{\|v\|^2}\om_0 - 
\frac{2}{\|v\|^3}dv \wedge \al \\
& = \om_0 + \pi^*\Om \\
& = \om_1.
\end{align*}
\end{proof}

\begin{remark}
The previous lemma shows that the symplectic manifold
$$ \bigg(\bigcup_{a \leq \mu \leq b}H^{-1}(\mu),\om_1\bigg) $$ 
has disconnected contact-type boundary for any $0<a<1/2<b$. It
can be compared with McDuff's example \cite{McD}. 
\end{remark}

It was proved in \cite{CGK,Lu1} that the tangent bundle of a compact manifold
endowed with any twisted symplectic form is geometrically bounded. Thus, consider the almost
complex structure $J$ on $TM$ such that $g(\cdot,\cdot)=\om_1(J\cdot,\cdot)$
defines a Riemannian metric with bounded sectional curvature and positive
injectivity radius.

By the previous lemma, given any $\ep>0$ sufficiently small, the subset $M_\ep
:= (TM\setminus M_0)\setminus \ov U_\ep$ has contact type boundary, that is,
its boundary is convex in the sense of Eliashberg-Gromov \cite{EG}. Thus, we
can consider the set $N_\ep := M_\ep \sqcup \partial M_\ep \times [0,\infty)$
(that we can identify with $TM\setminus M_0$) and a symplectic form $\rho_\ep$
on $N_\ep$ such that $\rho_\ep$ coincides with $\om_1$ on $M_\ep$ and the
vector field $\partial/\partial t$ is an extension of $X$ to $\partial M_\ep
\times [0,\infty)$ (that is, $-\partial/\partial t$ coincides with $-X$ on
$\partial M_\ep$) such that $L_{\partial/\partial t}\rho_\ep = \rho_\ep$. Note
that $\rho_\ep$ is exact, since $\om_1$ on $M_\ep$ is exact and $\rho_\ep$
restricted to $\partial M_\ep \times [0,\infty)$ is obviously exact.

Consider now the almost complex structure $J_\ep$ on $N_\ep$ given by
$J|_{M_\ep}$ on $M_\ep$ and by the pushforward $d\vr_t \circ J \circ d\vr_{-t}$
of $J$ by the flow $\vr$ of $\partial/\partial t$ on $\partial M_\ep \times
[0,\infty)$. By the construction, $J_\ep$ is invariant by $\vr$ on $\partial
M_\ep \times [0,\infty)$.

Define $g_\ep(\cdot,\cdot) = \rho_\ep(J_\ep\cdot,\cdot)$. By our previous
discussion, we have that $g_\ep$ defines a geometrically bounded Riemannian
metric on $M_\ep$. We need to prove that it also defines a geometrically
bounded Riemannian metric on $\partial M_\ep \times [0,\infty)$.

But, we have that $\vr_t^*\rho_\ep = e^t \rho_\ep$, that is,
$$ (\rho_\ep)_{(x,t)} = e^t (\vr_t)_*(\rho_\ep)_{(x,0)}.$$
Consequently, since $J_\ep$ is invariant by $\vr$ on $\partial M_\ep \times
[0,\infty)$, we obtain that 
$$(g_\ep)_{(x,t)} = e^t(\vr_t)_*(g_\ep)_{(x,0)}.$$
Thus, the sectional curvature $K_{(x,t)}(v,w)$ of $g_\ep$ on $(x,t)$ is given
by 
$$ e^{-t}K_{(x,0)}(d\vr_{-t}v,d\vr_{-t}w) $$ 
and the injectivity radius $i(N_\ep,(x,t))$ is greater than or equal to
$$ e^t i(N_\ep,(x,0)).$$
But, $g_\ep$ has bounded sectional curvature and positive injectivity radius on
$M_\ep$, proving the proposition.
\end{proof}

\begin{proof}[Proof of Theorem \ref{thmB}: ]
By our previous discussion, the magnetic flow restricted to an energy level
$k>1/2$ is topologically conjugated to the geodesic flow such that it has no
projected contractible periodic orbits. Consequently, we have that
$$ c_{\HZ}^G(U_k,\om_1) = \infty $$
for every $k>1/2$.

On the other hand, we have that for energies $k<1/2$ the (reparametrized) magnetic flow defines
a free Hamiltonian circle action (of period 1) on $U_k$ of infinite order on
$TM\setminus M_0$. By the previous proposition, given $\ep>0$ sufficiently
small, there exists a geometrically bounded symplectic form $\rho_\ep$ on $(TM\setminus M_0)$
that coincides with $\om_1$ on $(TM\setminus M_0)\setminus \ov U_\ep$. Thus, by
Theorem \ref{thmA} we have that for every $\delta>0$ there exists a finite constant $C>0$ such that, for every $\ep>0$ sufficiently small,
$$ c_{\HZ}^G(U_{1/2-\delta}\setminus \ov U_\ep,\om_1) = 
c_{\HZ}^G(U_{1/2-\delta}\setminus \ov U_\ep,\rho_\ep)
\leq C. $$
But,
$$ c_{\HZ}^G(U_{1/2-\delta},\om_1) = \sup_{0<\ep<1/2} c_{\HZ}^G(U_{1/2-\delta}\setminus \ov
U_\ep,\om_1) \leq C. $$
\end{proof}

\section{Proof of Theorem \ref{thmC}}

Before we give the proof of Theorem \ref{thmC}, we need some previous
lemmas. The first one can be derived by the arguments in the proof of
Proposition \ref{tame} and the fact $M\sm\Sig$ has a convex boundary.
However, for the reader's convenience we present an independent proof using the
structure of Stein manifolds:

\begin{lemma}
\label{tame1}
Given a neighborhood $U$ of $\Sig$, there exists an exact geometrically bounded symplectic form
$\rho$ on $M\sm\Sig$ such that $\rho$ coincides with $\Om$ on the complement of
$U$.
\end{lemma}

\begin{proof}
Note that $\Sig$ can be obtained as the zero set of a holomorphic section $s:
M \to L$ of a complex line bundle defined by the divisor $\Sig$. It can be show
that the function $\vr: M\sm\Sig \to \R$ given by
$$ \vr(x) = -\frac{1}{4\pi k}\log(\|s(x)\|^2) $$
is plurisubharmonic \cite{BC1}. The next lemma proved in \cite{BC2} shows that
we can construct a new plurisubharmonic function on $M\sm\Sig$ that coincides
with $\vr$ on $M\sm U$ and whose gradient vector field is complete.

\begin{lemma}[Lemma (3.1) of Biran-Cieliebak \cite{BC2}]
Let $(V,J,\vr)$ be a Stein manifold. Then for every $R \in \R$ there exists an
exhausting plurisubharmonic function $\vr_R: V \to \R$ with the following
properties:

1. $\vr_R = \vr$ on $V_{\vr\leq R} := \{x \in V; \vr(x)\leq R\}$;

2. the gradient vector field of $\vr_R$ is complete;

3. $\text{Crit}(\vr_R) = \text{Crit}(\vr)$ and for every $p \in
\text{Crit}(\vr_R)$, $\text{index}_p(\vr_R) = \text{index}_p(\vr)$.

In particular, the inclusion $(V_{\vr\leq R},\om_\vr) \hookrightarrow
(V,\om_{\vr_R})$ is a symplectic embedding.
\end{lemma}

Now, let $R>0$ be sufficiently great and consider the Riemannian metric
$g_{\vr_R}(\cdot,\cdot) = \om_{\vr_R}(J\cdot,\cdot)$, where $\om_{\vr_R} =
d(J^*d\vr_R)$.

Since $M$ is closed, we have that $g_{\vr_R}$ defines a geometrically bounded
Riemannian metric on $V_{\vr\leq R}$, that is, the sectional curvature of
$g_{\vr_R}$ restricted to $V_{\vr\leq R}$ is bounded from above and the
injectivity radius is positive. We need to prove that it also defines a
geometrically bounded Riemannian metric on all $M\sm\Sig$.

Given $x \in (M\sm\Sig)\sm V_{\vr\leq R}$, there exists $y \in \partial V_{\vr\leq
R}$ such that $X_{\vr_R}^t(y) = x$ for some $t>0$, where $X_{\vr_R}^t$ is the
gradient flow of $\vr_R$ with respect to $g_{\vr_R}$. But, we have that
$\vr_t^*\Om_{\vr_R} = e^t\Om_{\vr_R}$, that is,
$$ (\Om_{\vr_R})_{(x,t)} = e^t (\vr_t)_*(\Om_{\vr_R})_{(x,0)}.$$
Consequently, since the complex structure $J$ is invariant by $X_{\vr_R}^t$, we
obtain that 
$$ (g_{\vr_R})_{(x,t)} = e^t(\vr_t)_*(g_{\vr_R})_{(x,0)}.$$
Thus, the sectional curvature $K_{(x,t)}(v,w)$ of $g_{\vr_R}$ on $(x,t)$ is
given by 
$$ e^{-t}K_{(x,0)}(d\vr_{-t}v,d\vr_{-t}w) $$ 
and the injectivity radius $i(N_\ep,(x,t))$ is greater than or equal to
$$ e^t i(N_\ep,(x,0)).$$
But, $g_{\vr_R}$ has bounded sectional curvature and positive injectivity radius on
$V_{\vr\leq R}$, proving that $\Om_{\vr_R}$ is geometrically bounded.

Finally, define $\rho = \Om_{\vr_R}$ for $R>0$ sufficiently large.
\end{proof}

\begin{lemma}[Lemma (3.2) of Biran-Cieliebak \cite{BC2}]
\label{nullset}
Suppose that $M$ is subcritical and let $U \subset M\sm\Sig$ be an open subset
with compact closure. Let $\rho$ be a geometrically bounded symplectic form on $M\sm\Sig$ given
by the Lemma \ref{tame1}. Then there exists a symplectomorphism isotopic to
the identity $\psi: (M\sm\Sig,\rho) \to (M\sm\Sig,\rho)$ such that $\psi(U)
\cap \Delta = \emptyset$.
\end{lemma}

\begin{proof}
Remember that the isotropic CW-complex $\Delta$ is given by the union of the
stable manifolds of the gradient vector field $X_\vr$ of a certain
plurisubharmonic function $\vr$ with a complete gradient flow such that
$d(J^*d\vr) = \rho$ \cite{Bir1}.

Since $M$ is subcritical, we have that $\dim \Delta < \frac{1}{2}\dim_\R M$ and
so there exists a Hamiltonian isotopy $k_t: (M\sm\Sig,\rho) \to (M\sm\Sig,\rho)$
compactly supported in an arbitrarily small neighborhood of $\Delta$ such that
$k_1(\Delta) \cap \Delta = \emptyset$. Since $\Delta$ is compact, there exists a
neighborhood $V$ of $\Delta$ such that $k_1(V) \cap V = \emptyset$.

Fix $T$ so large that $X_\vr^{-T}(U) \subset V$. Since $k_1$ moves $V$ away from
itself, we have that
$$ X_\vr^T \circ k_1 \circ X_\vr^{-T}(U) \cap \Delta = \emptyset, $$
because $X_\vr^T$ leaves $(M\sm\Sig)\sm\Delta$ invariant. Since $X_\vr$ is
conformally symplectic, we have that
$$ \psi = X_\vr^T \circ k_1 \circ X_\vr^{-T} $$
is the desired symplectomorphism.
\end{proof}

\begin{proof}[Proof of Theorem \ref{thmC}:]
By Theorem \ref{Bir}, there exists an isotropic CW-complex
$\Delta \subset M$ (given by the union of the stale manifolds of the gradient
vector field associated to a plurisubharmonic function defined on
$M\setminus\Sigma$) whose complement $E:=M\setminus\Delta$ is symplectomorphic
to a standard symplectic disc bundle $(E_0,\frac{1}{k}\om_{\text{can}})$ over $\Sigma$
modeled on the normal bundle $N_\Sigma$ of $\Sigma$ in $M$, where
$$ \om_{\text{can}} = k\pi^*(\Om|_\Sigma) + d(r^2\al). $$
Moreover, this symplectomorphism $F:(E,\Om) \to (E_0,\frac{1}{k}\om_{\text{can}})$
sends $\Sigma$ to the zero section of $E_0$.

Note that $(E_0\sm\Sig,\frac{1}{k}\om_{\text{can}})$ has an obvious
Hamiltonian free circle action (of period $2\pi$) generated by the Hamiltonian
$H_0(x) = (1/2\pi k)\|x\|^2$ such that $H_1 = H_0 \circ F$ also generates a
Hamiltonian free circle action on $(E\sm\Sig,\Om)$.

By the Lefschetz theorem \cite{Bott}, we have that the inclusion $i:H_1^{-1}(r)
\to M\sm\Sig$ induces an injective map on the fundamental group such that the
order of the circle action on $E\sm\Sig \subset M\sm\Sig$ is equal to the order
of the circle action $X^{\om_{\text{can}}}_{H_0}$ generated by $H_0$ on $E_0\sm\Sig$.

Since the Chern class of $E_0$ is given by $k[\Om|_\Sig]$, we have that if
$[\Om]_{\pi_2(M)}=0$ then, by the Proposition \ref{propeulerclass} and the
Lefschetz theorem, the order of the action generated by $H_1$ is infinite.

On the other hand, the order of the action is always grater than or equal to
$k$. In effect, note that the $S^1$-bundle $S^1 \lra H_1^{-1}(r) \lra \Sig$
is equivalent to the quotient by $\Z_k$ of a $S^1$-bundle $S^1 \lra P \lra
\Sig$ with Euler class $[\Om_\Sig]$ such that the subgroup of
$\pi_1(H_1^{-1}(r))$ generated by the orbits of the action has order greater
than or equal to $k$.

Let $\rho_\ep$ be the exact geometrically bounded symplectic form on $M\sm\Sig$ given by the
Lemma \ref{tame1} that coincides with $\Om$ on the complement $M\sm U_\ep$ of
a tubular neighborhood of radius $\ep$ of $\Sig$. Note that
$(E\sm\Sig,\rho_\ep)$ also has a Hamiltonian free circle action that coincides
with $X^\Om_{H_1}$ on $M\sm U_\ep$.

Now, let $U \overset{i}{\hra} E\sm\Sig$ be an open subset with compact closure
and $\ep>0$ sufficiently small such that $U \subset E\sm U_\ep$. Then, by the
Theorem \ref{thmA},we have that if $[\Om]_{\pi_2(M)}=0$ then,
$$ c^{i^{-1}_*G}_{\HZ}(U,\Om) = c^{i^{-1}_*G}_{\HZ}(U,\rho_\ep) \leq 2\pi\|H_1|_U\| \leq
2\pi\|H_1|_E\| = 1/k. $$
where $G \subset \pi_1(E\sm\Sig) \simeq \pi_1(M\sm\Sig)$ is the subgroup
generated by the homotopy class of the orbits of $X^\Om_{H_1}$. On the other
hand, if $[\Om]_{\pi_2(M)} \neq 0$ and $k>2$ we have the inequality
$$ c^{i^{-1}_*G}_{\HZ}(U,\Om) = c^{i^{-1}_*G}_{\HZ}(U,\rho_\ep) \leq \frac{1}{k} + \frac{2}{k^2-2k}. $$

Finally, suppose that $(M,\Sig)$ is subcritical and let $U \overset{i}{\hra}
M\sm\Sig$ be an open subset with compact closure and $\ep>0$ sufficiently small
such that $U \subset M\sm U_\ep$. Then, by the Lemma \ref{nullset}, there
exists a symplectomorphism isotopic to the identity $\psi: (M\sm\Sig,\rho_\ep)
\to (M\sm\Sig,\rho_\ep)$ such that $\psi(U) \cap \Delta = \emptyset$. Take
$\ep>0$ sufficiently small such that $\psi(U) \subset M\sm U_\ep$. Thus, by the
weak monotonicity property, we conclude that
$$ c^{i^{-1}_*G}_{\HZ}(U,\Om) = c^{i^{-1}_*G}_{\HZ}(U,\rho_\ep) = 
c^{(\psi\circ i)^{-1}_*G}G_{\HZ}(\psi(U),\rho_\ep) \leq 1/k. $$ 
\end{proof}

\section{Proof of Theorem \ref{thmD}}

\begin{proof}[Proof of Theorem \ref{thmD}:]
The proof of Theorem \ref{thmD} is very similar to the proof of
Theorem \ref{thmC} but instead of complex hypersurfaces we consider symplectic
hypersurfaces of Donaldson type \cite{Don}.

In fact, by the symplectic neighborhood theorem, there exists a neighborhood
$V$ of $\Sig$ and a symplectomorphism $F:(V,\Om) \to (W,\om_{\text{can}})$ between
$(V,\Om)$ and a neighborhood $W$ of the zero section of the canonical disc
bundle $(E_0,\frac{1}{k}\om_{\text{can}})$, where
$$ \om_{\text{can}} = k\pi^*(\Om|_\Sigma) + d(r^2\al). $$
Moreover, $F$ sends $\Sigma$ to the zero section of $E_0$.

We will need the following lemma which generalizes the Lemma \ref{tame1} when
the complement of $\Sig$ may not be a Stein manifold:

\begin{lemma}
\label{tame2}
Given a neighborhood $U$ of $\Sig$, there exists a geometrically bounded symplectic form $\rho$
on $M\sm\Sig$ such that $\rho$ coincides with $\Om$ on the complement of $U$.
\end{lemma}

\begin{proof}
Note that
$$ (1/k)\om_{\text{can}} = \pi^*(\Om|_\Sigma) + (1/k)d(r^2\al) =
d(((1/k)r^2-1)\al). $$
Consequently, we have that the vector field
$$ X = \frac{(1/k)r^2-1}{(1/k)2r}\partial_r $$
is conformally symplectic with respect to $\om_{\text{can}}$.

Now, notice that $X$ point towards the interior of the compact domains bounded
by the hypersurfaces $r=\text{constant}$, that is, the subset
$E_\ep:=E_0\setminus\ov U_\ep$ has a convex boundary in the sense of
Eliashberg-Gromov \cite{EG}, where $U_\ep = \{x \in E_0;\ \ H_0(x) < \ep\}$ and
$H_0(x)$ is the radial coordinate of $x$.

Thus, we can consider the set $N_\ep := E_\ep \sqcup \partial E_\ep \times
[0,\infty)$  and a symplectic form $\rho_\ep$ on $N_\ep$ such that $\rho_\ep$
coincides with $\om_{\text{can}}$ on $E_\ep$ and the vector field $\partial_r$ is a
complete extension of $X$ to $\partial E_\ep \times [0,\infty)$ (that is,
$-\partial_r$ coincides with $-X$ on $\partial E_\ep$) such that
$L_{\partial_r}\rho_\ep = \rho_\ep$.

Now, we can proceed as in the proof of Proposition \ref{tame} to conclude
that $g_\ep(\cdot,\cdot) := \rho_\ep(J_\ep\cdot,\cdot)$ defines a geometrically
bounded Riemannian metric on $N_\ep$.
\end{proof}

In what follows we will repeat many arguments of the proof of Theorem
\ref{thmC}:

\vskip .5cm

Note that $(W\sm\Sig,\frac{1}{k}\om_{\text{can}})$ has an obvious Hamiltonian
free circle action (of period $2\pi$) generated by the Hamiltonian $H_0(x) =
(1/2\pi k)\|x\|^2$ such that $H_1 = H_0 \circ F$ also generates a Hamiltonian
free circle action on $(V\sm\Sig,\Om)$.

By the Lefschetz theorem for symplectic hypersurfaces of Donaldson type
\cite{Don}, there exists $k_0>0$ such that if $k>k_0$ then the inclusion
$i:H_1^{-1}(r) \to M\sm\Sig$ induces an injective map on the fundamental group
and so the order of the circle action on $V\sm\Sig \subset M\sm\Sig$ is equal
to the order of the circle action $X^{\om_{\text{can}}}_{H_0}$ generated by
$H_0$ on $W\sm\Sig$. Let $k_0 \geq 2$.

Since the Chern class of $E_0$ is given by $k[\Om|_\Sig]$, we have that, as in
the proof of Theorem \ref{thmC}, if $[\Om]_{\pi_2(M)}=0$ then the order
of the action generated by $H_1$ is infinite. On the other hand, the order of
the action is always grater than or equal to $k$.

Let $\rho_\ep$ be the geometrically bounded symplectic form on $M\sm\Sig$ given by the Lemma
\ref{tame2} that coincides with $\Om$ on the complement $M\sm U_\ep$ of a
tubular neighborhood of radius $\ep$ of $\Sig$.

Note that $(M\sm\Sig,\rho_\ep)$ is weakly exact, because $\rho_\ep$ is exact on
$V\sm\Sigma$ and, by the Lefschetz theorem, the inclusion $V\sm\Sig
\hookrightarrow M\sm\Sig$ induces a sobrejective application $\pi_2(V\sm\Sig)
\to \pi_2(M\sm\Sig)$.

Now, let $U \overset{i}{\hra} V\sm\Sig$ be a relative compact open subset and
$\ep>0$ sufficiently small such that $U \subset V\sm U_\ep$. Then, by the
Theorem \ref{thmA}, we have that if $[\Om]_{\pi_2(M)}=0$ then,
$$ c^{i^{-1}_*G}_{\HZ}(U,\Om) = c^{i^{-1}_*G}_{\HZ}(U,\rho_\ep) \leq 2\pi\|H_1|_U\| \leq
2\pi\|H_1|_E\| = 1/k. $$
where $G \subset \pi_1(V\sm\Sig)$ is the subgroup generated by the homotopy
class of the orbits of $X^\Om_{H_1}$. Analogously, if $[\Om]_{\pi_2(M)} \neq
0$ and $k>k_0 \geq 2$ we have the inequality 
$$ c^{i^{-1}_*G}_{\HZ}(U,\Om) = c^{i^{-1}_*G}_{\HZ}(U,\rho_\ep) \leq \frac{1}{k} + \frac{2}{k^2-2k}. $$
\end{proof}

\section{Proof of Theorem \ref{thmE}}
\label{proofthmE}

A diffeomorphism of a symplectic manifold $(M,\om)$ is called Hamiltonian if it
is given by the time one map of a time-dependent Hamiltonian vector field. Let
$Ham(M,\om)$ be the group of Hamiltonian diffeomorphisms of $(M,\om)$. It was
discovered by Hofer $\cite{Hof}$ that this group carries a natural biinvariant
Finsler metric with a non-degenerate distance function. This distance is
defined as follows. The Lie algebra $\cal A$ of $Ham(M,\om)$ consists of all
smooth functions on $M$ satisfying a normalization condition which ensures
that, when $M$ is open, $F \in \cal A$ iff $F$ has compact support, and, when
$M$ is closed, $F \in \cal A$ iff $F$ has zero mean with respect to the volume
form induced by $\om$. The adjoint action of $Ham(M,\om)$ on $\cal A$ is the
standard action of diffeomorphisms on functions. Consider the $L^\infty$-norm
in $\cal A$,
$$ \|F\| = \max_M F - \min_M F $$
This norm is invariant under the adjoint action, and thus defines a biinvariant
Finsler metric on $Ham(M,\om)$. Now, let $\{\ga_t\}$, $t \in [0,1]$, a path of
Hamiltonian diffeomorphisms, that is, a curve $\ga: [0,1] \to Ham(M,\om)$. Let
$F(x,t)$ be its normalized Hamiltonian function, that is, $F(\cdot,t) \in \cal
A$ for all $t$. We define then the length of this curve by
$$ l(\ga) = \int_0^1 \|F(\cdot,t)\| dt $$
The Hofer's metric is then defined by
$$ \rho(\phi,\psi) = \inf_\ga l(\ga) $$
where the infimum is taken over all paths $\ga$ which join $\phi$ and $\psi$.

For a subset $U$ of $M$ define the displacement energy $e(U)$ as the measure of
the distance between the identity map and the set of $\psi \in Ham(M,\om)$ which
displaces $U$ from itself, in the sense that $\psi(U) \cap U \not= \emptyset$.

Using the geometry of the Hofer's metric and that, for a manifold whose Euler
characteristic vanishes, the displacement energy of the zero section of $T^*M$,
with the twisted symplectic form given by a non-vanishing magnetic field, is
equal to zero, Polterovich \cite{Pol1} proved the following theorem:

\begin{theorem} [Polterovich \cite{Pol1}]
Let $M$ be a closed manifold whose Euler characteristic vanishes and $\Om$ a
non-vanishing weakly exact magnetic field. Then there exist contractible
magnetic closed orbits on a sequence of arbitrarily small energy levels.
\end{theorem}

The idea of the proof is the following \cite{Pol1}: take a smooth function
$r(x)$, $x \in [0,\infty)$, which equals $x-2\ep$ on $[0,\ep]$, vanishes on
$[3\ep,\infty)$ and is strictly increasing on $[0,3\ep)$. Thus, the closed
orbits of the (normalized) Hamiltonian $F=r \circ H$ corresponds to
reparametrized magnetic periodic orbits whose energy is less than $3\ep$. The
key point is that the minimum set of $F$ coincides with the zero section of
$T^*M$ and so its displacement energy vanishes \cite{Pol2}. It implies that the
asymptotic non-minimality  of $F$ is strictly less than 1 and thus that $F$ (as
a vector in the Lie algebra $\cal A$) does not generate a minimal geodesic in
$Ham(M,\om)$ \cite{Pol1}. Then, using a result of Lalonde and McDuff
\cite{LMcD}, we conclude that the Hamiltonian $F$ has a non-constant
contractible closed orbit, since $\Om$ is weakly exact.

\begin{remark}
It should be noted that the contractibility of the periodic orbits given by
Lalonde and McDuff follows by the fact that the proof of Theorem 5.4 of
\cite{LMcD} relies on the Hofer-Zehnder capacity-area inequality (Theorem 1.17
of \cite{LMcD}) proved by Floer, Hofer and Viterbo \cite{FHV}, where the
periodic orbits are contractible (see \cite{McDSl}).
\end{remark}

\begin{proof}[Proof of Theorem \ref{thmE}:]
Consider the trivial $S^1$-bundle $S^1 \lra P = M \times S^1
\overset{\pi}{\lra} M$, the flat connection $d\theta$ and the product metric
$\bar g = \lg,\rg$. Let $\Om_P$ be the pullback of the magnetic field $\Om$ to
$P$. The lifted circle action to $TP$ is Hamiltonian with respect to the
twisted symplectic form $\om_1 = \om_o + \pi_P^*\Om_P$, where $\pi_P:TP \to P$
is the canonical projection. In fact, let $\wh{X}$ the vector field on $TP$
that generates the circle action. Then,
\begin{align*} i_{\wh{X}}\om_1(\xi) &
=i_{\wh{X}}\om_0(\xi)+i_{\wh{X}}\pi_P^*\pi^*\Om(\xi) \\ & =
i_{\wh{X}}\om_0(\xi) + \Om(\pi_*X,\pi_*(\pi_P)_*\xi) \\ & =
i_{\wh{X}}\om_0(\xi) = d\theta(\xi),
\end{align*}
where $d\theta$ above is viewed as a function in $TP$. Moreover, the form
$d\theta$ is invariant by the magnetic flow of $TP$, because taking a magnetic
orbit $\ga$ in $P$, we have that 
\begin{align*} \ga^\pr d\theta(\ga^\pr) & = \ga^\pr\lg\ga^\pr,X\rg \\ & =
\lg\ga^{\pr\pr},X\rg + \lg\ga^\pr,\nabla_{\ga^\pr}X\rg \\ & = \lg
Y_P(\ga^\pr),X \rg \\ & = \Om_P(\ga^\pr,X) \\ & = \Om(\pi_*\ga^\pr,\pi_*X) = 0,
\end{align*} 
where $Y_P$ denotes the Lorentz force of the magnetic field $\Om_P$ (defined as
the bundle map $Y:TP \to TP$ uniquely characterized by the property that
$(\Om_P)_x(v,w)= \bar g_x(Y_x(v),w)\ \forall v,w \in T_xP$) and in the third
equation we used the fact that $X$ is parallel.

Thus, the magnetic flow on $TP$ is equivariant with respect to the Hamiltonian
lifted action of $S^1$ to $TP$. On the other hand, it is easy to see that the
Marsden-Weinstein reduced symplectic manifold with respect to the circle action
on $(TP,\om_1)$ is given by $(TM,\om_0 + \pi_M^*\Om)$.

Moreover, the pullback of the Hamiltonian metric $H_g$ in $TM$ by
$d\pi|_{d\theta^{-1}(\xi)}$ is given by $H_{\bar g}|_{d\theta^{-1}(\xi)} -
\frac{1}{2}|\xi|^2$. Consequently, the reduced dynamics of the magnetic flow
$\phi_P$ of $P$ on $d\theta^{-1}(\xi)$ is given by the magnetic flow $\phi_M$
on $M$.

Now, suppose that $\Om$ is weakly exact. Since $P$ is the product, the fibers
of $P$ are non-contractible and, in fact, the subgroup generated by the fibers
in $\pi_1(P)$ is isomorphic to $\Z$. Since the Euler characteristic of $P$
vanishes and $\pi^*\Om$ is weakly exact (because $\pi_2(M \times S^1)
\simeq \pi_2(M)$), we can apply Polterovich's theorem to conclude that
$\phi_P$ has contractible closed orbits $\ga_n$ of arbitrarily small energy.
Since the subgroup in $\pi_1(P)$ generated by the fibers of $P$ is cyclic
infinite, these orbits cannot be given by the fibers. This ensures that they
are projected on non-constant contractible periodic orbits of $\phi_M$.
Finally, note that the energy of these projected orbits is equal to
$$ H_{\bar g}(\ga_n) - \frac{1}{2}|\xi|_n^2 \overset{n\to\infty}{\lra} 0 $$
where $\xi_n = d\theta(\ga_n^\pr) \overset{n\to\infty}{\lra} 0$, because
$\|\ga_n^\pr\|_{\bar g} \to 0$ as $n$ grows to infinite.
\end{proof}


\begin{thebibliography}{ABCD}

\bibitem{Ar1} V. Arnold, {\em Some remarks on flows of line elements and
frames}, Soviet Math. Dokl. {\bf 2} (1961), 565--564.

\bibitem{Ar2} V. Arnold, {\em First steps of symplectic topology}, Russian Math.
Surveys {\bf 41} (1986), 1--21.

\bibitem{Bir1} P. Biran, {\em Lagrangian barriers and symplectic embeddings},
GAFA {\bf 11} (2001), no. 3, 407--464.

\bibitem{BC1} P. Biran, K. Cieliebak, {\em Symplectic topology on subcritical
manifolds}, Comment. Math. Helv. {\bf 76} (2001), 712--753.

\bibitem{BC2} P. Biran, K. Cieliebak, {\em Lagrangian embeddings into subcritical
Stein manifolds}, Israel Journal of Mathematics {\bf 127} (2002), 221--244.

\bibitem{BPS} P. Biran, L. Polterovich, D. Salamon, {\em Propagation in
Hamiltonian dynamics and relative symplectic homology}, preprint 2001,
math.SG/0108134.

\bibitem{Bott} R. Bott, {\em On a theorem of Lefschetz}, Michigan Math. J. {\bf
6} (1959), 211--216.

\bibitem{CGK} K. Cieliebak, V. Ginzburg, E. Kerman. {\em Symplectic homology and
periodic orbits near symplectic submanifolds}, preprint 2002 math.DG/0210468.

\bibitem{CMP} G. Contreras, L. Macarini, G. Paternain, {\em Periodic orbits for
exact magnetic flows on surfaces}, to appear in Int. Math. Res. Notices.

\bibitem{Don} S. Donaldson, {\em Symplectic submanifolds and almost-complex
geometry}, J. Differential Geom. {\bf 44} (1996), no. 4, 666--705.

\bibitem{EG} Y. Eliashberg, M. Gromov, {\em Convex symplectic manifolds}
Several complex variables and complex geometry, Part 2 (Santa Cruz, CA, 1989), 
135-162, Proc. Sympos. Pure Math., 52, Part 2, Amer. Math. Soc., Providence,
RI, 1991.

\bibitem{EGH} Y. Eliashberg, A. Givental, H. Hofer, {\em Introduction to
symplectic field theory}, GAFA - Special Volume, Part II (2000), 560--673.

\bibitem{FHV} A. Floer, H. Hofer, C. Viterbo, {\em The Weinstein conjecture in
$P \times \C^l$}, Math. Z. {\bf 203} (1990), 469--482.

\bibitem{Gin} V. Ginzburg, {\em A charge in a magnetic field: Arnold's problems
1981-9, 1982-24, 1984-4, 1994-14, 1996-17, and 1996-18}, preprint (disponible at http://count.ucsc.edu/~ginzburg).

\bibitem{Gin2} V. Ginzburg, {\em The Hamiltonian Seifert conjecture: examples
and open problems}, math.DG/0004020.

\bibitem{GK} V. Ginzburg, E. Kerman, {\em Periodic orbits in magnetic fields in
dimensions greater than two}, in {\em Geometry and topology in dynamics},
Contemp. Math {\bf 246}, 119--121.

\bibitem{Gro} M. Gromov, {\em Pseudo holomorphic curves in symplectic
manifolds}, Inv. Math. {\bf 82} (1985), 307--347.

\bibitem{Her} M. Herman, {\em Exemples de flots Hamiltoniens dont aucune
perturbation en topologie $C^\infty$ n'a d'orbites periodiques sur un ouvert de
surfaces d'energies}, Comptes Rendus de l'Acad\e mie des Sciences de Paris {\bf
412} (1991), 989--994.

\bibitem{Hof} H. Hofer, {\em Estimates for the energy of a symplectic map}, Com.
Math. Helvetici {\bf 68} (1993), 48--72.

\bibitem{HV} H. Hofer, C. Viterbo, {\em The Weinstein conjecture in the
presence of holomorphic spheres}, Comm. Pure. Apl. Math. {\bf 45} (1992),
583--622.

\bibitem{HV2} H. Hofer, C. Viterbo, {\em The Weinstein conjecture in cotangent
bundles and related results}, Ann. Sc. Norm. Sup. Pisa {\bf 15} (1988), 411--445.

\bibitem{HZ} H. Hofer, E. Zehnder, {\em Symplectic invariants and Hamiltonian
dynamics}, Birkh\"auser Advanced Texts; Basel-Boston-Berlin, 1994.

\bibitem{Ji} M. Jiang, {\em Hofer-Zehnder symplectic capacity for 2-dimensional
manifolds}, Proc. Royal Soc. Edinburgh {\bf 123A} (1993), 945--950.

\bibitem{Ker} E. Kerman, {\em Symplectic geometry and the motion of a particle
in a magnetic field}, PhD. thesis, University of California - Santa Cruz, 2000.

\bibitem{Kob} S. Kobayashi, {\em Principal fibre bundles with the 1-dimensional
toroidal group}, Tohoku J. Math. {\bf 8} (1956), 29-45.

\bibitem{LMcD} L. Lalonde, D. McDuff, {\em Hofer's $L^\infty$ geometry:
geodesics and stability I, II}, Inv. Math. {\bf 122} (1995), 1--33,35--69.

\bibitem{Lu1} G. Lu, {\em The Weinstein conjecture on some symplectic manifolds
containing the holomorphic spheres}, Kyushu J. Math {\bf 52} (1998), 331--351.

\bibitem{Lu2} G. Lu, {\em Errata to the "The Weinstein conjecture on some
symplectic manifolds containing the holomorphic spheres"}, Kyushu J. Math {\bf
54} (2000).

\bibitem{Lu3} G. Lu, {\em A note on Lagrangian barrier theorem by P. Biran}, math.SG/0111183.

\bibitem{Mac}  L. Macarini, {\em Hofer-Zehnder semicapacity of cotangent bundles and symplectic submanifolds}, math.SG/0303230.

\bibitem{MS} L. Macarini, F. Schlenk, {\em A refinement of the Hofer--Zehnder theorem on the existence 
of closed trajectories near a hypersurface}, math.SG/0305148.

\bibitem{McD} D. McDuff, {\em Symplectic manifolds with contact-type
boundaries}, Inventiones Mathematicae {\bf 103} (1991), no. 3, 651--671.

\bibitem{McDSl} D. McDuff, J. Slimowitz, {\em Hofer-Zehnder capacity and length
minimizing Hamiltonian paths}, Geom. Topol. {\bf 5} (2001), 799--830 

\bibitem{Nov} S. Novikov, {\em The Hamiltonian formalism and many-valued
analogue of Morse theory} Russian Math.surveys {\bf 37} (1982), 1--56.

\bibitem{PS} R. Palais, T. Stewart, {\em Torus bundles over a torus}, Proc. of
AMS {\bf 42} (1961), 29-29.

\bibitem{Pol1} L. Polterovich, {\em Geometry on the group of Hamiltonian
diffeomorphisms}, Proceedings of the ICM - 1998 (Vol. II), 401--410.

\bibitem{Pol2} L. Polterovich, {\em An obstacle to non-Lagrangian
intersections}, pp. 575-586 in The Floer Memorial Volume edited by H. Hofer, C.
Taubes, A. Weinstein and E. Zehnder, Birkh\"auser, 1995.

\bibitem{Sch} M. Schwarz, {\em On the action spectrum for closed symplectically
aspherical manifolds},  Pacific J. Math. {\bf 193} (2000), no. 2, 419--461.

\bibitem{Si} K. Siburg, {\em Symplectic capacities in two dimensions},
Manuscripta Math. {\bf 78} (1993), 149--163.

\bibitem{Str} M. Struwe, {\em Existence of periodic solutions of Hamiltonian
systems on almost every energy surface}, Bol. da Sociedade Brasileira de
Matem\'atica {\bf 20} (1990), 49--58.

\bibitem{Vit} C. Viterbo, {\em Exact Lagrange submanifolds, periodic orbits and
the cohomology of free loop spaces},  J. Differential Geom. {\bf 47} (1997),
no. 3, 420--468.

\bibitem{Vit2} C. Viterbo, {\em Functors and computations in Floer homology
with applications. I.}, Geom. Funct. Anal. {\bf 9} (1999),  no. 5, 985--1033. 

\bibitem{We} A. Weinstein, {\em On the hypothesis of Rabinowitz's periodic
orbit theorems}, J. Diff. Equation {\bf 33} (1979), 353--358.

\bibitem{Zeh} E. Zehnder, {\em Remarks on periodic solutions on hypersurfaces},
NATO ASI Series, Series C, In: Periodic Solutions of Hamiltonian Systems and
Related Topics, {\bf 209} (1987), 267--279.

\end{thebibliography}
\end{document}